\documentclass[a4paper,11pt,fleqn]{article}
\usepackage{amsmath,amssymb,amsthm,bm,mathtools}
\usepackage{geometry}
\usepackage[hidelinks]{hyperref}
\usepackage{enumitem,booktabs,array}
\usepackage[most]{tcolorbox}
\usepackage{algorithm}
\usepackage{algpseudocode}
\usepackage{float}
\usepackage{graphicx}

\setlist[enumerate]{topsep=3pt,itemsep=3pt,parsep=0pt}
\allowdisplaybreaks[2]
\newcommand{\R}{\mathbb{R}}
\newcommand{\norm}[1]{\left\| #1\right\|}
\newcommand{\ip}[2]{\left\langle #1,#2\right\rangle}
\newcommand{\argmin}{\operatorname*{arg\,min}}

\newcommand{\prox}{\operatorname{prox}}
\newcommand{\J}{\nabla F}

\theoremstyle{plain}
\newtheorem{theorem}{Theorem}[section]
\newtheorem{lemma}[theorem]{Lemma}
\newtheorem{proposition}[theorem]{Proposition}
\newtheorem{corollary}[theorem]{Corollary}
\theoremstyle{definition}

\newtheorem{remark}[theorem]{Remark}

\newtcolorbox{mybox}[1]{
  title={#1},
  colback=white,
  colframe=black,
  colbacktitle=white,
  coltitle=black,
  fonttitle=\bf,
  sharp corners,
  boxrule=0.5pt,
  left=8pt,
  right=8pt,
  top=8pt,
  bottom=8pt
}

\begin{document}
\title{Adaptive Linesearch-Free Prox-Linear Method for Convex Composite Optimization}
\author{Kotaro Matsuda \thanks{Graduate School of Informatics, Kyoto University, Japan, \texttt{matsuda.kotaro.28e@st.kyoto-u.ac.jp}}
        \and
        Nobuo Yamashita \thanks{Graduate School of Informatics, Kyoto University, Japan, \texttt{nobuo@i.kyoto-u.ac.jp}}}
\date{November 27, 2026}
\maketitle

\begin{abstract}

We consider the convex composite optimization problem
$\min_{x\in\R^n} g(F(x))$,
where $g:\R^m\to\R$ is a finite-valued,
possibly nonsmooth convex function that is
nondecreasing in each coordinate,
and $F:\R^n\to\R^m$ has continuously
differentiable convex components.
We propose an adaptive linesearch-free prox-linear method,
denoted by $\mathrm{AdGDC}(a,b)$, that determines its stepsizes
using information obtained from successive iterations,
without requiring prior knowledge of Lipschitz constants
or performing line searches.
Building on adaptive gradient methods for smooth convex optimization,
the proposed method incorporates two nonnegative quantities
arising from convexity into a unified stepsize update rule,
with their contributions controlled by the parameters $(a,b)$.
We establish a fundamental inequality and prove boundedness
of the generated iterates under basic convexity assumptions.
Under the additional assumption that the gradients of the
inner functions are locally Lipschitz continuous,
we prove a uniform positive lower bound on the stepsizes,
an $O(1/K)$ convergence rate for the best iterate,
convergence of the objective values to the optimal value,
and convergence of the entire sequence to an optimal solution.
Numerical experiments on nonsmooth convex composite problems
based on the maximum of log-cosh loss functions demonstrate
the effectiveness of the proposed method compared with
existing prox-linear methods.

\end{abstract}


\section{Introduction}
\label{sec:introduction}

We consider the convex composite optimization problem
\begin{equation}\label{eq:problem_intro}
 \min_{x\in\R^n} h(x):=g(F(x)),
\end{equation}
where each component $F_i$ of $F=(F_1,\ldots,F_m):\R^n\to\R^m$ is a convex function of class
$C^1$, and $g:\R^m\to\R$ is a finite-valued, possibly nondifferentiable convex function that is
nondecreasing in each coordinate. 
Under these assumptions, $h$ is continuous, finite-valued, and convex. 
We further assume that problem~\eqref{eq:problem_intro} admits an optimal solution.
We use the following notation for the optimal value and the solution set:
\begin{equation}\label{eq:solution_set}
  h^* :=\min_{x\in\R^n} h(x)>-\infty,\qquad  X^* :=\argmin_{x\in\R^n} h(x)\ne\emptyset.
\end{equation}

The class of problems of the form~\eqref{eq:problem_intro} includes smooth unconstrained convex
optimization and accommodates both nonsmooth regularization and the aggregation of multiple
objectives. 
For example, setting $m=1$ and $g(v)=v$ gives the usual smooth convex optimization problem
$\min_{x\in\R^n} F_1(x)$.
For a smooth convex function $f$ and $\lambda\ge0$, the $\ell_1$-regularized problem
\[
 \min_{x\in\R^n}\bigl\{f(x)+\lambda\norm{x}_1\bigr\}
\]
also admits this formulation.
Indeed, it suffices to set
\[
 F(x)=(f(x),x_1,-x_1,\ldots,x_n,-x_n),\qquad
 g(u)=u_1+\lambda\sum_{i=1}^n
 \max\{u_{2i},u_{2i+1}\}.
\]
Choosing $g(u)=\max_{1\le i\le m}u_i$ gives a min--max problem in which one minimizes the maximum
of several nonlinear convex functions.
This maximum-type formulation also encompasses certain instances of the Pascoletti--Serafini (PS)
model for the multiobjective optimization problem $\min F(x)$.
For example, given a reference point $c\in\R^m$ and a positive direction vector $r\in\R^m$, $r>0$,
setting
\[
 g(u)=\max_{1\le i\le m}\frac{u_i-c_i}{r_i}
\]
yields the scalarized objective of the PS model~\cite{BanholzerMechelliVolkwein2022}.

In problem~\eqref{eq:problem_intro}, even when the inner functions $F_i$ are differentiable, the
objective $h=g\circ F$ is generally nonsmooth.
Consequently, standard gradient descent cannot be applied directly to $h$.
One way to address this nonsmoothness is to form a first-order approximation of the inner mapping
$F$ and solve a subproblem that retains the structure of the outer function $g$.
Since the first-order approximation of $F$ introduces an error, one must account simultaneously
for the nonsmoothness of $g$ and the linearization error of $F$.
Efficient solution of such problems therefore calls for algorithms that exploit their composite
structure.

Various iterative methods for nonlinear composite optimization compute search directions using an
approximate model of the inner mapping.
These methods can be broadly classified according to whether they use only first derivatives of
the inner mapping or also use second-order information.
Methods of the former type linearize $F$ at the current point $x$ to obtain the model
$g\bigl(F(x)+\J(x)^\top d\bigr)$.
Minimizing this model without regularization gives the subproblem
\begin{equation}
  \min_{d\in\R^n} g\bigl(F(x)+\J(x)^\top d\bigr)
  \label{eq:subproblem_intro}.
\end{equation}
Gauss--Newton-type methods are representative of approaches based on this subproblem. 
Burke and Ferris~\cite{BurkeFerris1995} analyzed a Gauss--Newton method for convex-composite
optimization, establishing globalization by line search and local quadratic convergence under
suitable regularity conditions. 
Ochs and Malitsky~\cite{OchsMalitsky2019} proposed a model-function-based conditional gradient
method that minimizes the model over a constraint set and also treated Gauss--Newton-type methods
for nonlinear composite problems.
For the unregularized subproblem~\eqref{eq:subproblem_intro}, however, existence and uniqueness of
a solution need not hold.
This has motivated the prox-linear method, which adds a positive quadratic regularization term to
the linearized model and minimizes the resulting subproblem
objective~\cite{DrusvyatskiyPaquette2019,LewisWright2016}.
Models incorporating second-order information have also been proposed.
For example, Burke and Engle~\cite{BurkeEngle2020} analyzed local convergence of Newton and
quasi-Newton methods for composite optimization problems involving the composition of a convex
piecewise linear-quadratic function with a smooth mapping.
Such methods can capture the shape of the function more accurately than a first-order
approximation, although computing or approximating second derivatives and solving the resulting
subproblems can be costly.


At the current point $x$, the prox-linear method considered here adds a proximal regularization
term to the linearized model of the inner mapping and solves the subproblem
\[
 \min_{d\in\R^n}
 \left\{
 g\bigl(F(x)+\J(x)^\top d\bigr)
 +\frac{1}{2t}\norm{d}^2
 \right\}.
\]
Its solution $s$ is used to update $x$ to $x+s$. 
Here $\J(x)$ is the Jacobian of $F$, and $t>0$ is the stepsize. 
The subproblem uses a first-order approximation of the inner mapping $F$ and retains the outer
function $g$ in its original form. 
The quadratic proximal term controls the length of the search direction and makes the subproblem
strongly convex, ensuring uniqueness of its solution. 
As above, combining a first-order approximation of the inner mapping with proximal regularization
is the basic idea of the prox-linear method. 
Since the regularization coefficient $1/t$ controls the length of the search direction, the choice
of stepsize $t$ plays an important role in the performance of the algorithm.

Work on the prox-linear method and its extensions has examined convergence, complexity, model
accuracy, and adaptive parameter selection.
Lewis and Wright~\cite{LewisWright2016} studied a composite optimization framework with proximal
regularization, and Drusvyatskiy and Paquette~
\cite{DrusvyatskiyPaquette2019}
analyzed the complexity of the prox-linear method for a more general class of composite
minimization problems.
Bolte et al.~\cite{BolteChenPauwels2020} proposed the Multiprox method, which uses multiple
proximal terms reflecting the curvature of individual components, and established a convergence
rate for convex composite problems.
Mukkamala et al.~\cite{MukkamalaFadiliOchs2022} proposed a method based on the Model Approximation
Property, which bounds the model approximation error in terms of a Bregman distance.

Existing prox-linear methods commonly use a fixed stepsize based on a Lipschitz constant or a
related bound, or adjust the stepsize adaptively by line search.
In the former approach, the relevant constants may be difficult to estimate in advance, and a
global bound may yield an overly conservative stepsize.
In the latter, each iteration may require several function evaluations or subproblem solves.
To address these issues, Hu et al.~\cite{HuLiWangYangZhu2023} assumed that the outer convex
function attains its minimum and proposed a prox-linear method with an adaptive stepsize based on
this minimum value.
Yang et al.~\cite{YangShenHuLi2026} further developed its convergence analysis and globalization
procedures.
Although adaptive methods for nonlinear composite problems have thus received some attention,
there remains scope for developing and analyzing update rules that use information from the
iterates to determine the stepsize without a prior Lipschitz estimate or a line search.

Our aim is to propose a convergent prox-linear method for problem~\eqref{eq:problem_intro} that
chooses its stepsize adaptively without a line search or prior knowledge of a Lipschitz constant
for the gradient.
Our starting point is the work of Malitsky and Mishchenko~
\cite{MalitskyMishchenko2020,MalitskyMishchenko2024}
on Adaptive gradient descent without descent (AdGD) for smooth convex optimization and its
subsequent extensions.
AdGD updates the stepsize explicitly using gradient differences between iterates, and its
convergence is guaranteed under local Lipschitz continuity of the gradient.
Latafat et al.~
\cite{LatafatThemelisPatrinos2024,LatafatThemelisStellaPatrinos2025}
proposed the Adaptive proximal gradient method (AdaPG), which extends the stepsize update using,
among other quantities, an inner product involving the gradient difference.
Various other extensions of adaptive methods have been studied~
\cite{OikonomidisEtAl2024,ZhouMaYang2026,PandeyZhangZheng2026}.
These studies focus primarily on smooth objectives or on additive composite problems whose
objectives are sums of smooth and nonsmooth functions.
In the nonlinear composite problems considered here, the outer function $g$ may be nonsmooth, and
the linearization error of the inner mapping $F$ must also be taken into account.
Thus, these update rules do not apply directly.

We propose the Adaptive gradient descent method for composite problems (AdGDC), which extends
existing adaptive gradient methods by incorporating properties specific to the composite structure
into stepsize selection.
The proposed method, AdGDC$(a,b)$, determines the stepsize explicitly from information obtained
during the iteration.
The parameters $(a,b)$ combine two nonnegative quantities arising from convexity, providing a
common treatment of update rules that use different nonnegative terms.
In particular, when $g$ is affine, suitable parameter choices recover the stepsize updates of
AdGD~\cite[Algorithm~1]{MalitskyMishchenko2024} and AdaPG~\cite{LatafatThemelisPatrinos2024}.
We prove boundedness of the iterates under basic convexity assumptions.
Assuming in addition that each $\nabla F_i$ is locally Lipschitz continuous, we establish a
uniform positive lower bound on the stepsizes, an $O(1/K)$ best-iterate convergence rate,
convergence of the objective values, and convergence of the entire sequence of iterates.

The paper is organized as follows. 
Section~\ref{sec:related_research} reviews the adaptive gradient methods for smooth problems most
closely related to our work. 
Section~\ref{sec:algorithm_AdGDC} defines the subproblem and stepsize update for AdGDC$(a,b)$.
Sections~\ref{sec:gen_ineq} through~\ref{sec:xk_conv} establish, in turn, the fundamental
inequality and boundedness of the iterates, the stepsize lower bound and convergence rate, and
convergence of the entire sequence.
Section~\ref{sec:parameter_variants} discusses representative parameter choices.
Section~\ref{sec:experiments} presents numerical experiments, and Section~\ref{sec:conclusion}
concludes.

We conclude this section by fixing our notation.
Let $\R^n$ denote real $n$-dimensional space and $\R^n_+$ its nonnegative orthant.
For $u,v\in\R^n$, let $\ip{u}{v}$ denote the standard inner product and $\norm{u}$ the Euclidean
norm.
For a matrix $A\in\R^{n\times m}$, let $\norm{A}$ denote the operator norm induced by the
Euclidean norm.
For $x\in\R^n$ and $r\ge0$, let
$\overline{B}(x,r):=\{y\in\R^n:\norm{y-x}\le r\}$
denote the closed Euclidean ball centered at $x$
with radius $r$.
For $c\in\R$, define $[c]_+:=\max\{c,0\}$.
We denote the Jacobian of $F=(F_1,\ldots,F_m):\R^n\to\R^m$ by
\[
 \J(x):=
 [\nabla F_1(x),\ldots,\nabla F_m(x)]
 \in\R^{n\times m}.
\]
We use superscripts for vectors and subscripts for scalars indexed by the iteration number $k$.

\section{Related work on Adaptive Gradient Methods}
\label{sec:related_research}
For smooth convex optimization problems, 
Malitsky and Mishchenko~\cite{MalitskyMishchenko2020,MalitskyMishchenko2024} proposed AdGD, a linesearch-free
method that determines its stepsize from gradient differences.
Latafat et al.~\cite{LatafatThemelisPatrinos2024} extended the AdGD stepsize update to obtain
AdaPG, which permits larger stepsizes.
We review the results of Malitsky and Mishchenko~\cite{MalitskyMishchenko2024} and Latafat et
al.~\cite{LatafatThemelisPatrinos2024} most closely related to our method and its analysis.

\subsection{AdGD of Malitsky and Mishchenko}
\label{sec:malitskyandmishcenko}
Consider the smooth convex optimization problem
\begin{equation}\label{eq:prior_smooth_problem}
 \min_{x\in\R^n} f(x),
\end{equation}
where $f$ is convex and differentiable and its solution set is nonempty.
For the gradient iteration applied to $f$,
\[x^{k+1}:=x^k-t_k\nabla f(x^k),\]
define, for $k\ge1$,
\[
 s^k:=x^k-x^{k-1},\qquad
 y^k:=\nabla f(x^k)-\nabla f(x^{k-1}),\qquad
 \theta_k:=\frac{t_k}{t_{k-1}}.
\]
Malitsky and Mishchenko~\cite[Algorithm~1]{MalitskyMishchenko2024} considered AdGD with an
arbitrary initial stepsize $t_0>0$ and the following update for $k\ge1$:
\begin{equation}\label{eq:prior_mm2024_alg1}
 t_k:=\min\left\{
 \sqrt{1+\theta_{k-1}}\,t_{k-1},\quad
 \frac{\norm{s^k}}{\sqrt{2}\norm{y^k}}
 \right\}.
\end{equation}
When $\norm{y^k}=0$, the second term is interpreted as $+\infty$.
In~\cite[Algorithm~1]{MalitskyMishchenko2024}, the authors set $\theta_0=0$.

For comparison with our subsequent analysis, let
\begin{align*}
  x^*\in\argmin f,\quad
  \delta_k:=f(x^k)-f(x^*), \quad
  \sigma_k^f:=f(x^{k-1})-f(x^k)+\ip{\nabla f(x^k)}{s^k}\ge0.
\end{align*}
Convexity of $f$ gives $\ip{y^k}{s^k}\ge0$.

The second term in~\eqref{eq:prior_mm2024_alg1} ensures that
\begin{equation}\label{eq:prior_mm_condition}
 2t_k^2\norm{y^k}^2\le\norm{s^k}^2
\end{equation}
holds.
Using this condition and convexity, the authors of~\cite{MalitskyMishchenko2024} showed that
for every $x^*\in \argmin f$,
\begin{equation}
 \norm{x^{k+1}-x^*}^2+\norm{s^{k+1}}^2
 +2t_k(1+\theta_k)\delta_k \le
 \norm{x^k-x^*}^2+\norm{s^k}^2
 +2t_k\theta_k\delta_{k-1}
 \label{eq:prior_mm_onestep}
\end{equation}
holds.
This is equation~(13) of~\cite{MalitskyMishchenko2024} with $\gamma=1/\sqrt{2}$.

These estimates yield the following convergence result.

\begin{theorem}[Convergence of AdGD]
\label{thm:prior_mm_convergence}
Let $f:\R^n\to\R$ be convex, let $\nabla f$ be locally Lipschitz continuous, and suppose that
$\argmin f\ne\emptyset$.
For arbitrary $x^0\in\R^n$ and $t_0>0$, consider AdGD with $\theta_0=0$ and stepsizes updated
by~\eqref{eq:prior_mm2024_alg1}.
Unless the method terminates at an optimal solution after finitely many iterations, the generated
sequence $\{x^k\}$ converges to some $x^* \in\argmin f$.

Moreover, for any $x^*\in\argmin f$, define
\[
 R_{x^*}^2:=
 \norm{x^0-x^*}^2+
 2t_0^2\norm{\nabla f(x^0)}^2.
\]
Then $\{x^k\}$ lies in $\overline{B}(x^*,R_{x^*})$.
Furthermore, there exists $t_{\inf}>0$ such that $t_k\ge t_{\inf}$ for all $k$, and for every
$K\ge1$ we have
\[
 \min_{1\le k\le K}\bigl(f(x^k)-f(x^*)\bigr)
 \le
 \frac{R_{x^*}^2}{2\sum_{k=1}^K t_k}
 \le
 \frac{R_{x^*}^2}{2t_{\inf}K}
 =O\left(\frac{1}{K}\right).
\]
\end{theorem}

The preceding theorem follows by an argument analogous to the convergence proof for~
\cite[Algorithm~2]{MalitskyMishchenko2024}
and is described as a consequence of that argument.
The authors do not state it separately for~
\cite[Algorithm~1]{MalitskyMishchenko2024}
in their paper.

To clarify the connection with our analysis, we now show how to refine
the estimate~\eqref{eq:prior_mm_onestep}.
First, the gradient update and convexity give
\begin{align}
 &\norm{x^{k+1}-x^*}^2+\norm{s^{k+1}}^2
 +2t_k(1+\theta_k)\delta_k
 +2t_k\theta_k
 \bigl(\ip{y^k}{s^k}+\sigma_k^f\bigr)\notag\\
 &\qquad\le
 \norm{x^k-x^*}^2+2t_k\theta_k\delta_{k-1}
 +2t_k^2\norm{y^k}^2
 \label{eq:prior_mm_refined}.
\end{align}
We call this inequality, and its counterpart in the analysis of the stepsize rule proposed in
Section~\ref{sec:algorithm_AdGDC}, the fundamental inequality.
Section~\ref{sec:meaning_stepsize_ineq} gives the detailed derivation
of~\eqref{eq:prior_mm_refined} and its consequence~\eqref{eq:prior_mm_onestep}.

Using~\eqref{eq:prior_mm_condition} and $\sigma_k^f\ge0$ further yields
\begin{align}
 &\norm{x^{k+1}-x^*}^2+\norm{s^{k+1}}^2
 +2t_k(1+\theta_k)\delta_k
 +2t_k\theta_k\ip{y^k}{s^k}\notag\\
 &\qquad\le
 \norm{x^k-x^*}^2+\norm{s^k}^2
 +2t_k\theta_k\delta_{k-1}
 \label{eq:prior_mm_strong}.
\end{align}
Since $\ip{y^k}{s^k}\ge0$,
\eqref{eq:prior_mm_strong} gives
a sharper estimate than~\eqref{eq:prior_mm_onestep}.

\subsection{AdaPG of Latafat et al.}
Latafat et al.~\cite{LatafatThemelisPatrinos2024} proposed the adaptive proximal gradient method
AdaPG$^{q,r}$, which uses an inner product involving the gradient difference. 
In particular, applying the method with $(q,r)=(1,1/2)$ to the smooth
problem~\eqref{eq:prior_smooth_problem} gives the following update that uses
the nonnegative term $\ip{y^k}{s^k}$ retained in~\eqref{eq:prior_mm_strong}
to enlarge the stepsize:
\[
 t_k:=\min\left\{
 \sqrt{1+\theta_{k-1}}\,t_{k-1},\quad
 \frac{\norm{s^k}}{\sqrt{2\left[
 \norm{y^k}^2-\ip{y^k}{s^k}/t_{k-1}
 \right]_+}}
 \right\}.
\]
When the denominator vanishes, the second term is interpreted as $+\infty$.
Since $\ip{y^k}{s^k}\ge0$, this rule
permits a larger stepsize upper bound than~\eqref{eq:prior_mm2024_alg1}.
For more general convex additive composite problems, Latafat et
al.~\cite[Theorem~1.1]{LatafatThemelisPatrinos2024} established convergence properties analogous
to those of Theorem~\ref{thm:prior_mm_convergence} under assumptions including local Lipschitz
continuity of $\nabla f$.
\medskip

The preceding studies primarily concern $\min f(x)$ or the problem of minimizing the sum of a
smooth convex function $f_{\mathrm{s}}$ and a nonsmooth convex function $f_\mathrm{n}$, namely
$\min\{f_{\mathrm{s}}(x)+f_{\mathrm{n}}(x)\}$.
In the composite problem~\eqref{eq:problem_intro} considered here, $\min\; g(F(x))$, we linearize
the nonlinear mapping $F$ and retain the outer convex function $g$ when solving the subproblem.
We must therefore account simultaneously for the nonsmoothness of $g$ and the linearization error
of $F$.
We propose AdGDC as an extension of the AdGD-type update to this problem structure.

\section{The proposed method}
\label{sec:algorithm_AdGDC}

In this section, we linearize $F$ in the objective $g(F(x))$ to define a subproblem and propose
the algorithm AdGDC$(a,b)$ based on this subproblem.
Section~\ref{sec:subproblem} defines the linearized model of the objective and regularizes it to
obtain the subproblem used in the algorithm.
Section~\ref{sec:stepsize_rule} details the stepsize selection rule.
Section~\ref{sec:algorithm_AdGDC_all} presents the complete algorithm.

\subsection{Prox-linear subproblem}
\label{sec:subproblem}

For the convex composite optimization problem~\eqref{eq:problem_intro}, define the linearized
model at $x$ by
\begin{equation}\label{eq:model}
 m_x(d):=g\bigl(F(x)+\J(x)^\top d\bigr).
\end{equation}
Convexity of the functions $F_i$ and coordinatewise monotonicity of $g$ imply that, for all
$x,y\in\R^n$,
\begin{equation}\label{eq:model_minorant}
 m_x(y-x)=g(F(x)+\nabla F(x)^\top (y-x))\le g(F(y)) = h(y).
\end{equation}
Moreover, by definition,
\[
  m_x(0)=h(x) \qquad(x\in\R^n).
\]

Since $g$ is finite-valued and convex on $\R^m$, $\partial g(v)$ is nonempty and its subgradients
are uniformly bounded on every bounded set.
Coordinatewise monotonicity of $g$ also gives
\begin{equation}\label{eq:nonnegative_z}
 \partial g(v)\subset\R_+^m\qquad(v\in\R^m).
\end{equation}

We next define the subproblem used in the algorithm.
Given $k\ge0$ and a stepsize $t_k>0$, consider
\begin{equation}
 s^{k+1}:=\argmin_{d\in\R^n}
 \left\{m_{x^k}(d)+\frac{1}{2t_k}\norm d^2\right\}
 \label{eq:subproblem}.
\end{equation}
Using its solution $s^{k+1}$, we update the iterate by setting
\[
  x^{k+1}:=x^k+s^{k+1}.
\]

If the subproblem~\eqref{eq:subproblem} has an optimal solution, the chain rule for finite-valued
convex functions gives the optimality condition
\begin{equation}
  0\in \partial m_{x^k}(s^{k+1}) +\frac{s^{k+1}}{t_k} = \J(x^k)\partial g(\widetilde{F}^k)+\frac{s^{k+1}}{t_k} 
  \label{eq:optim_condition},
\end{equation}
where $\widetilde{F}^k$ is defined by
\begin{equation}\label{eq:uk}
 \widetilde{F}^k:=F(x^k)+\J(x^k)^\top s^{k+1}.
\end{equation}

The following theorem establishes existence of a subproblem solution and provides a stopping
criterion.
\begin{theorem}\label{thm:subproblem}
For each $t_k>0$, problem~\eqref{eq:subproblem} has a unique solution.
Moreover, there exists $z^k\in\partial g(\widetilde{F}^k)\subseteq \R_+^m$ such that
\begin{equation}\label{eq:dual_opt}
 \J(x^k)z^k=-s^{k+1}/t_k.
\end{equation}
If $s^{k+1}=0$, then $x^k\in X^* $.
\end{theorem}
\begin{proof}
Since $m_{x^k}$ is finite-valued and convex, adding the positive quadratic term gives a strongly
convex subproblem with a unique minimizer.
Its optimality condition is
\[ 0\in\J(x^k)\partial g(\widetilde{F}^k)+s^{k+1}/t_k,\] 
which implies~
\eqref{eq:dual_opt}.
If $s^{k+1}=0$,
then~\eqref{eq:model_minorant} and the subgradient inequality give, for every $y$,
\[ h(y)\ge h(x^k)+\ip{\J(x^k)z^k}{y-x^k}=h(x^k), \] 
hence $x^k\in X^* $.
\end{proof}

By Theorem~\ref{thm:subproblem}, one may, for example, stop the algorithm in practice when
$\|s^k\|$ is sufficiently small.

\subsection{Stepsize update rule}
\label{sec:stepsize_rule}

Motivated by the stepsize~\eqref{eq:prior_mm2024_alg1} of Malitsky and Mishchenko, we propose an
explicit rule for updating the stepsize $t_k$ in subproblem~\eqref{eq:subproblem} for
problem~\eqref{eq:problem_intro}.
The first term $\sqrt{1+\theta_{k-1}}\,t_{k-1}$ in the minimum is the same as
in~\eqref{eq:prior_mm2024_alg1}.
We modify the second term to account for the composite structure of~\eqref{eq:problem_intro}.

For $k\ge1$, fix the subgradient $z^{k-1}\in \partial g (\widetilde{F}^{k-1})$ appearing in the
optimality condition of the preceding iteration and define
\begin{align}
 y^k&:=\bigl(\J(x^k)-\J(x^{k-1})\bigr)z^{k-1},\label{eq:yk}\\
 \eta_k&:=g(F(x^k))-g(\widetilde{F}^{k-1})
      -\ip{z^{k-1}}{F(x^k)-\widetilde{F}^{k-1}}\ge0,\label{eq:eta_k} \\
\sigma_k&:=
 h(x^{k-1})-h(x^k)+\eta_k
 +\ip{\J(x^k)z^{k-1}}{s^k}\ge0\label{eq:sigma_k}.
\end{align}

The quantity $y^k$ is the difference of the Jacobians weighted by the fixed subgradient $z^{k-1}$
from the preceding iteration.
The quantity $\eta_k$ is the error, at $F(x^k)$, of the first-order approximation of $g$ at
$\widetilde{F}^{k-1}$ formed using a subgradient.
In the differentiable case it corresponds to the Bregman distance.
Since $z^{k-1}\in\partial g(\widetilde{F}^{k-1})$, the subgradient inequality gives $\eta_k\ge0$.
Moreover, $\sigma_k$ is the error at $x^{k-1}$ of the first-order approximation of the objective
$h$ at $x^k$ using $\J(x^k)z^{k-1}$ and $\eta_k$.
We prove $\sigma_k \ge 0$ in Section~\ref{sec:gen_inec_derivation},
Lemma~\ref{lem:approx_subgrad}.

\begin{remark}
  \label{rem:eta0_g_affine}
Since $g$ need not be differentiable, $y^k$ is in general different from
$\nabla h(x^k)-\nabla h(x^{k-1})$.
If, for example, $g$ is affine, then $\eta_k=0$.
Moreover, $y^k=\nabla h(x^k)-\nabla h(x^{k-1})$.
\end{remark}

By Theorem~\ref{thm:subproblem}, we have $z^{k-1}\ge0$.
Together with the convexity of each $F_i$, this implies
\begin{equation}\label{eq:ymonotone}
 \ip{y^k}{s^k}
 =\sum_{i=1}^m z_i^{k-1}
 \ip{\nabla F_i(x^k)-\nabla F_i(x^{k-1})}{x^k - x^{k-1}}\ge0.
\end{equation}

We now combine $\ip{y^k}{s^k}$ and $\sigma_k$ in selecting the stepsize.
Define the set of admissible parameters $a, b$ by
\begin{equation}
  \Omega:=\{(a,b)\in[0,1]\times[0,1] \colon a+b<2\}
  \label{eq:def_Omega},
\end{equation}
and fix any $(a,b)\in\Omega$.
We denote the gradient method using the stepsize determined by the pair $(a,b)$ by AdGDC$(a,b)$.

For $(a,b)\in \Omega$ and $k\ge1$, define
\begin{equation}
  \rho_k(a,b):=a\,\ip{y^k}{s^k}+b\,\sigma_k \ge 0
  \label{eq:def_rho}.
\end{equation}
The inequality $\rho_k(a,b) \ge 0$ follows from~
\eqref{eq:ymonotone} and the bound $\sigma_k \ge 0$ established in Lemma~\ref{lem:approx_subgrad}.
Although $\rho_k (a,b)$ depends on $(a,b)$, the parameters $(a,b)$ are fixed for a given
algorithm.
We therefore write simply $\rho_k$ whenever $(a,b)$ is clear from the context.

Given $s^k$, $y^k$, $\eta_k$, and $\rho_k$, we update $\widehat t_k$ and $t_k$ as follows.
Take an arbitrary initial value $t_0>0$ and set $\theta_0=1$.

\begin{mybox}{$\widehat t_k$ in AdGDC$(a,b)$}
\vspace{-10pt}
\begin{equation}\label{eq:root}
 \widehat t_k:=
 \displaystyle\frac{\norm{s^k}^2}{2\eta_k+\sqrt{4\eta_k^2+
       2\left[\norm{y^k}^2+\dfrac{\eta_k-\rho_k}{t_{k-1}}\right]_+
       \norm{s^k}^2}}.
\end{equation}
If the denominator is $0$, set $\widehat t_k:=+\infty$.
\end{mybox}
\vspace*{-14pt}
\begin{equation}\label{eq:stepsize}
 t_k:=\min\left\{\sqrt{1+\theta_{k-1}}\,t_{k-1},\ \widehat t_k\right\},
 \qquad\theta_k:=t_k/t_{k-1}.
\end{equation}

The derivation of the update~\eqref{eq:root} for $\widehat t_k$ begins in
Section~\ref{sec:gen_ineq}.
The form of the update~\eqref{eq:stepsize} is the same as in the method of Malitsky and
Mishchenko.
The new element proposed here is~\eqref{eq:root}.

\subsection{The complete algorithm}
\label{sec:algorithm_AdGDC_all}

Algorithm~\ref{alg:adgdc} gives the complete AdGDC$(a,b)$ method.

\begin{algorithm}[H]
\caption{AdGDC$(a,b)$}
\label{alg:adgdc}
\begin{algorithmic}[1]
\Require Parameters $(a,b)\in\Omega$, initial point $x^0\in\R^n$, initial stepsize $t_0>0$,
initial $\theta$ given by $\theta_0\gets 1$

  \For{$k=0,1,2,\ldots$}

    \If{$k\ge 1$}
      \State
        $y^k\gets
        \bigl(\J(x^k)-\J(x^{k-1})\bigr)z^{k-1}$

      \State
        $\eta_k\gets
        g(F(x^k))-g(\widetilde{F}^{k-1})
        -\ip{z^{k-1}}{F(x^k)-\widetilde{F}^{k-1}}$

      \State
        $\sigma_k\gets
        h(x^{k-1})-h(x^k)+\eta_k
        +\ip{\J(x^k)z^{k-1}}{s^k}$

      \State
        $\rho_k\gets
        a\ip{y^k}{s^k}+b\sigma_k$

\State $\widehat t_k\gets$ the value given by~\eqref{eq:root}

      \State
        $t_k\gets
        \min\bigl\{
        \sqrt{1+\theta_{k-1}}\,t_{k-1},
        \ \widehat t_k
        \bigr\}$

      \State $\theta_k\gets t_k/t_{k-1}$
    \EndIf

    \State
      $\displaystyle s^{k+1}\gets
      \argmin_{d\in\R^n}
      \left\{
      g\bigl(F(x^k)+\J(x^k)^\top d\bigr)
      +\frac{\norm d^2}{2t_k}
      \right\}$

    \If{$s^{k+1}=0$}
      \State \Return $x^k$
    \EndIf

    \State $x^{k+1}\gets x^k+s^{k+1}$

    \State
      $\widetilde{F}^k\gets
      F(x^k)+\J(x^k)^\top s^{k+1}$

    \State
Choose $z^k\in\partial g(\widetilde{F}^k) \subseteq \R^m_+$ such that
      $\J(x^k)z^k=-s^{k+1}/t_k$

  \EndFor
\end{algorithmic}
\end{algorithm}

\begin{remark}
If $F$ is affine, AdGDC$(a,b)$ reduces to a
proximal point method for any $(a,b)\in\Omega$.
Indeed, since $F(x^k)+\J(x^k)^\top d=F(x^k+d)$,
we have $m_{x^k}(d)=h(x^k+d)$.
Consequently, the iteration takes the form
\[
  x^{k+1}=\prox_{t_kh}(x^k).
\]
\end{remark}

The following propositions relate AdGDC$(a,b)$ to
AdGD~\cite{MalitskyMishchenko2024} and AdaPG~\cite{LatafatThemelisPatrinos2024}.

\begin{proposition}[Relation to AdGD]\label{prop:special}
Suppose that the outer function $g$ is affine, that is,
$g(u)=\ip{\alpha}{u}+\beta$ for some
$\alpha\in\R_+^m$ and $\beta\in\R$.
Then AdGDC$(0,0)$ has the same stepsize update rule
as AdGD of Malitsky and Mishchenko~
\cite{MalitskyMishchenko2024}.
\end{proposition}

\begin{proof}
Let $(a,b)=(0,0)$.
When $g$ is affine, we have $z^{k-1}=\alpha$, $\eta_k=0$, and
$y^k=\nabla h(x^k)-\nabla h(x^{k-1})$.
Thus $\widehat t_k=\norm{x^k - x^{k-1}}/(\sqrt{2}\norm{\nabla h(x^k)-\nabla h(x^{k-1})})$ if
$y^k\ne0$, and $\widehat t_k=+\infty$ if $y^k=0$.
Since $s^{k+1}=-t_k\nabla h(x^k)$, this gives the same stepsize update as AdGD.
\end{proof}

\begin{proposition}[Relation to AdaPG]\label{prop:adapg_special}
Suppose that the outer function $g$ is affine, that is, $g(u)=\ip{\alpha}{u}+\beta$ for some
$\alpha\in\R_+^m$ and $\beta\in\R$.
Then, AdGDC$(1,0)$ has the same stepsize update as AdaPG of Latafat et
al.~\cite{LatafatThemelisPatrinos2024} with $(q,r)=(1,1/2)$.
\end{proposition}
\begin{proof}
Let $(a,b)=(1,0)$.
If $g$ is affine, then $z^{k-1}=\alpha$, $\eta_k=0$, and $y^k=\nabla h(x^k)-\nabla h(x^{k-1})$.
Consequently,
\eqref{eq:root} becomes
\[
 \widehat t_k =
 \frac{\norm{s^k}}{\sqrt{2\left[\norm{y^k}^2-\ip{y^k}{s^k}/t_{k-1}\right]_+}}
\]
with the value $+\infty$ when the denominator is $0$.
This gives the same stepsize update as AdaPG with $(q,r)=(1,1/2)$.
The subproblem also reduces to the usual gradient step $s^{k+1}=-t_k\nabla h(x^k)$.
\end{proof}

Thus AdGDC extends AdGD of Malitsky and Mishchenko~\cite{MalitskyMishchenko2024} in the sense of
Proposition~\ref{prop:special} and AdaPG of Latafat et al.~\cite{LatafatThemelisPatrinos2024} in
the sense of Proposition~\ref{prop:adapg_special}.
For $b>0$, the stepsize is new even for the smooth minimization problem $\min f(x)$.

\section{The fundamental inequality and boundedness of the iterates}
\label{sec:gen_ineq}

For AdGDC$(a,b)$, we derive a fundamental inequality corresponding to the
estimate~\eqref{eq:prior_mm_refined} for smooth convex optimization presented in
Section~\ref{sec:related_research} and use it to establish boundedness of the iterates.
In doing so, we explain the origin of the stepsize rules~\eqref{eq:root} and~\eqref{eq:stepsize}.

Define the optimality gap by
\begin{equation}\label{eq:delta_def}
 \delta_k:=h(x^k)- h^* \ge0.
\end{equation}

\subsection{Derivation of the fundamental inequality}
\label{sec:gen_inec_derivation}

We derive the fundamental inequality, the counterpart of inequality~\eqref{eq:prior_mm_refined}
for AdGD of Malitsky and Mishchenko~\cite{MalitskyMishchenko2024}.
This inequality underlies our analysis of the proposed algorithm.
We first establish several auxiliary estimates.
Before proving them, we record the relations used in this section.
For every $k\ge1$, the optimality condition from the preceding iteration,
\eqref{eq:dual_opt}, yields
\begin{equation}
 t_{k-1}\J(x^{k-1})z^{k-1}=-s^k \notag.
\end{equation}
Substituting this into the definition~\eqref{eq:yk} of $y^k$ gives
\begin{equation}\label{eq:tw}
 t_k\J(x^k)z^{k-1}
 =-\theta_k s^k+t_ky^k.
\end{equation}

Lemma~\ref{lem:approx_subgrad} uses the subgradient obtained at the preceding iteration to bound
the objective and its linearized model from below.

\begin{lemma}\label{lem:approx_subgrad}
Let $k\ge1$.
For every $d\in\R^n$ and $w\in\R^n$,
\begin{align}
 m_{x^k}(d)&\ge h(x^k)+\ip{\J(x^k)z^{k-1}}{d}-\eta_k,
 \label{eq:approx_model}\\
 h(w)&\ge h(x^k)+\ip{\J(x^k)z^{k-1}}{w-x^k}-\eta_k
 \label{eq:approx_global}.
\end{align}
In particular, for $\sigma_k$ defined in~\eqref{eq:sigma_k},
\begin{equation}\label{eq:prev_gap}
 \sigma_k\ge0.
\end{equation}
\end{lemma}

\begin{proof}
By $z^{k-1}\in\partial g(\widetilde{F}^{k-1})$, for every $v\in\R^m$ we have
\[
 g(v)\ge g(\widetilde{F}^{k-1})
 +\ip{z^{k-1}}{v-\widetilde{F}^{k-1}}.
\]
Taking $v=F(x^k)+\J(x^k)^\top d$ and using the definition~\eqref{eq:model} of $m_x (d)$ gives
\begin{align*}
 m_{x^k}(d)
 &\ge g(\widetilde{F}^{k-1})
  +\ip{z^{k-1}}{
    F(x^k)+\J(x^k)^\top d-\widetilde{F}^{k-1}}\\
 &=g(\widetilde{F}^{k-1})
  +\ip{z^{k-1}}{F(x^k)-\widetilde{F}^{k-1}}
  +\ip{\J(x^k)z^{k-1}}{d}\\
 &=h(x^k)-\eta_k+\ip{\J(x^k)z^{k-1}}{d}.
\end{align*}
The last equality uses the definition~\eqref{eq:eta_k} of $\eta_k$.
This proves~\eqref{eq:approx_model}.

Next, by~\eqref{eq:model_minorant}, for every $w\in\R^n$ we have
\begin{align*}
 h(w)
 &\ge m_{x^k}(w-x^k)\\
 &\ge h(x^k)
   +\ip{\J(x^k)z^{k-1}}{w-x^k}-\eta_k,
\end{align*}
which gives~\eqref{eq:approx_global}.

Furthermore, substituting $w=x^{k-1}$ into~\eqref{eq:approx_global} and using $x^{k-1}-x^k=-s^k$
yields
\[
 h(x^{k-1})
 \ge h(x^k)-\ip{\J(x^k)z^{k-1}}{s^k}-\eta_k.
\]
Since $\sigma_k=h(x^{k-1})-h(x^k)+\eta_k+\ip{\J(x^k)z^{k-1}}{s^k}$ by definition,
\eqref{eq:prev_gap} follows.
\end{proof}

We briefly interpret Lemma~\ref{lem:approx_subgrad}.
For $\varepsilon\ge0$, define the $\varepsilon$-subdifferential of a convex function $h$ by
\[
 \partial_\varepsilon h(x):=
 \left\{p\in\R^n:\ h(w)\ge h(x)+\ip{p}{w-x}-\varepsilon
 \quad(\forall w\in\R^n)\right\}.
\]
Then~\eqref{eq:approx_global} means that
\[
 \J(x^k)z^{k-1}\in\partial_{\eta_k}h(x^k).
\]
In the smooth setting of Malitsky and Mishchenko~\cite{MalitskyMishchenko2020}, we have $\eta_k=0$
and this vector equals $\nabla h(x^k)$.
In our analysis of the prox-linear method, we must retain the error $\eta_k$ arising from the
nonlinear composite structure.

For use in Lemma~\ref{lem:distance} below, define $e^k\; (k\ge 1)$ by
\begin{equation}\label{eq:ek}
 e^k:=s^{k+1}+t_k\J(x^k)z^{k-1}
 =s^{k+1}-\theta_ks^k+t_ky^k.
\end{equation}
The quantity $e^k$ is the residual in the optimality condition~\eqref{eq:optim_condition} for the
current subproblem~\eqref{eq:subproblem} when $z^k$ is replaced by the preceding subgradient
$z^{k-1}$.

\begin{remark}
\label{rem:e0_g_affine}
If $g$ is affine, then $z^k=z^{k-1}$, so the subproblem optimality condition gives
\[
 e^k
 =s^{k+1}+t_k\J(x^k)z^{k-1}
 =t_k\J(x^k)(z^{k-1}-z^k)
 =0.
\]
\end{remark}

The following lemma relates $\|x^{k+1}-x^*\|^2$ to $\|x^{k}-x^*\|^2$ and will be used to derive
the fundamental inequality.
Here $x^* \in X^*$.

\begin{lemma}\label{lem:distance}
For every $k\ge1$ and $x^*\in X^* $,
\begin{equation}\label{eq:distance_residual}
 \norm{x^{k+1}-x^*}^2+2t_k\delta_k+\norm{e^k}^2
 \le\norm{x^k-x^*}^2
 +t_k^2\norm{\J(x^k)z^{k-1}}^2+2t_k\eta_k.
\end{equation}
\end{lemma}
\begin{proof}
The optimality condition~\eqref{eq:optim_condition} for subproblem~\eqref{eq:subproblem} and
convexity of $m_{x^k}$ give
\[
 m_{x^k}(x^*-x^k)
 \ge m_{x^k}(s^{k+1})
 -\frac1{t_k}\ip{s^{k+1}}{x^*-x^{k+1}}.
\]
By~\eqref{eq:model_minorant}, the left-hand side is at most $h(x^*)= h^* $.
Hence
\[
 \ip{s^{k+1}}{x^{k+1}-x^*}
 \le t_k\bigl( h^* -m_{x^k}(s^{k+1})\bigr).
\]
The identity
\[
 \norm{x^{k+1}-x^*}^2
 =\norm{x^k-x^*}^2
 +2\ip{s^{k+1}}{x^{k+1}-x^*}
 -\norm{s^{k+1}}^2
\]
then yields
\begin{equation}\label{eq:distance_intermediate}
 \norm{x^{k+1}-x^*}^2
 \le\norm{x^k-x^*}^2-\norm{s^{k+1}}^2
 +2t_k\bigl( h^* -m_{x^k}(s^{k+1})\bigr).
\end{equation}

Substituting $d=s^{k+1}$ into~\eqref{eq:approx_model} further gives
\[
  h^* -m_{x^k}(s^{k+1})
 \le-\delta_k
 -\ip{\J(x^k)z^{k-1}}{s^{k+1}}+\eta_k.
\]
Substitution into~\eqref{eq:distance_intermediate} yields
\[
 \norm{x^{k+1}-x^*}^2+2t_k\delta_k
 \le\norm{x^k-x^*}^2-\norm{s^{k+1}}^2-2t_k\ip{\J(x^k)z^{k-1}}{s^{k+1}}
 +2t_k\eta_k.
\]
On the right-hand side, \eqref{eq:ek} gives
\[
 -\norm{s^{k+1}}^2
 -2t_k\ip{\J(x^k)z^{k-1}}{s^{k+1}}=
 -\norm{e^k}^2
 +t_k^2\norm{\J(x^k)z^{k-1}}^2,
\]
and~\eqref{eq:distance_residual} follows.
\end{proof}

We now combine the preceding estimates to obtain the fundamental inequality for AdGDC.

\begin{lemma}[Fundamental inequality]\label{lem:malitsky}
For every $k\ge1$ and $x^*\in X^* $,
\begin{align}
 &\norm{x^{k+1}-x^*}^2+\norm{s^{k+1}}^2
 +2t_k(1+\theta_k)\delta_k
 +2\norm{e^k}^2
 +2t_k\theta_k \left(\ip{y^k}{s^k} + \sigma_k \right)
 \notag\\
 &\qquad\le\norm{x^k-x^*}^2+2t_k\theta_k\delta_{k-1}
 +2t_k^2\norm{y^k}^2+2t_k(2+\theta_k)\eta_k.
 \label{eq:refined_twostep}
\end{align}
\end{lemma}
\begin{proof}
First, the optimality condition~\eqref{eq:optim_condition} for subproblem~\eqref{eq:subproblem}
gives $-s^{k+1}/t_k\in\partial m_{x^k}(s^{k+1})$, so
\[
 h(x^k)=m_{x^k}(0)
 \ge m_{x^k}(s^{k+1})+\frac{\norm{s^{k+1}}^2}{t_k}.
\]
Combining this with~\eqref{eq:approx_model} yields
\[
 \norm{s^{k+1}}^2
 +t_k\ip{\J(x^k)z^{k-1}}{s^{k+1}}
 \le t_k\eta_k.
\]
On the other hand, expanding the square of~\eqref{eq:ek} gives
\[
 \norm{s^{k+1}}^2+\norm{e^k}^2
 =2\norm{s^{k+1}}^2
 +2t_k\ip{\J(x^k)z^{k-1}}{s^{k+1}}
 +t_k^2\norm{\J(x^k)z^{k-1}}^2.
\]
Hence
\begin{equation}\label{eq:residual_square}
 \norm{s^{k+1}}^2+\norm{e^k}^2
 \le t_k^2\norm{\J(x^k)z^{k-1}}^2+2t_k\eta_k.
\end{equation}

Next, \eqref{eq:tw} implies $\J(x^k)z^{k-1}=-s^k/t_{k-1}+y^k$.
Therefore,
\begin{equation}\label{eq:p_expand}
 t_k^2\norm{\J(x^k)z^{k-1}}^2
 =t_k^2\norm{y^k}^2+\theta_k^2\norm{s^k}^2
 -2t_k\theta_k\ip{y^k}{s^k}.
\end{equation}
Substituting the same expression into the definition of~\eqref{eq:sigma_k} gives
\[
 \frac{\norm{s^k}^2}{t_{k-1}}
 =\delta_{k-1}-\delta_k+\eta_k
 +\ip{y^k}{s^k}-\sigma_k.
\]
Multiplication of both sides by $t_k\theta_k$ yields
\[
 \theta_k^2\norm{s^k}^2
 =t_k\theta_k
 \bigl(\delta_{k-1}-\delta_k+\eta_k
 +\ip{y^k}{s^k}-\sigma_k\bigr).
\]
Substituting this into~\eqref{eq:p_expand} gives
\begin{equation}\label{eq:refined_curvature}
 t_k^2\norm{\J(x^k)z^{k-1}}^2
 =t_k^2\norm{y^k}^2
 +t_k\theta_k(\delta_{k-1}-\delta_k+\eta_k)-t_k\theta_k\ip{y^k}{s^k}
 -t_k\theta_k\sigma_k.
\end{equation}

Finally, adding the estimate~\eqref{eq:distance_residual} in Lemma~\ref{lem:distance}
to~\eqref{eq:residual_square} gives
\[
 \norm{x^{k+1}-x^*}^2+\norm{s^{k+1}}^2
 +2t_k\delta_k+2\norm{e^k}^2 \le\norm{x^k-x^*}^2
 +2t_k^2\norm{\J(x^k)z^{k-1}}^2+4t_k\eta_k.
\]
Substituting into the term $t_k^2\norm{\J(x^k)z^{k-1}}^2$ on the right-hand side
the expression in~\eqref{eq:refined_curvature} and rearranging gives~
\eqref{eq:refined_twostep}.
\end{proof}

\subsection{Stepsize selection and boundedness of the iterates}
\label{sec:meaning_stepsize_ineq}

When $g$ is affine, Remark~\ref{rem:eta0_g_affine} gives $\eta_k=0$ and
Remark~\ref{rem:e0_g_affine} gives $e^k=0$.
Thus inequality~\eqref{eq:refined_twostep} in Lemma~\ref{lem:malitsky} becomes
\begin{align*}
 &\norm{x^{k+1}-x^*}^2+\norm{s^{k+1}}^2
 +2t_k(1+\theta_k)\delta_k
 +2t_k\theta_k\left(\ip{y^k}{s^k}+\sigma_k\right)\\
 &\qquad\le
 \norm{x^k-x^*}^2+2t_k\theta_k\delta_{k-1}
 +2t_k^2\norm{y^k}^2.
\end{align*}
This corresponds to the estimate in Section~\ref{sec:malitskyandmishcenko},
namely~\eqref{eq:prior_mm_refined}.
Discarding the fourth term on the left-hand side,
$2t_k\theta_k(\ip{y^k}{s^k}+\sigma_k)$,
and imposing the stepsize condition
$2t_k^2\norm{y^k}^2\le\norm{s^k}^2$
used by Malitsky and Mishchenko~
\cite{MalitskyMishchenko2024}
yield~\eqref{eq:prior_mm_onestep}.

For the composite problem~\eqref{eq:problem_intro}, we may have $\eta_k>0$.
Stepsize selection must therefore account for the term involving $\eta_k$. 
We also use the fourth term on the left-hand side discarded by Malitsky and Mishchenko to allow
greater flexibility in the choice of stepsize.
Specifically, we propose choosing $t_k$ to satisfy
\begin{equation}\label{eq:curvature_condition}
 2t_k^2\norm{y^k}^2
 +2t_k\left(2+\frac{t_k}{t_{k-1}}\right)\eta_k
 -\frac{2t_k^2}{t_{k-1}}\rho_k
 \le\norm{s^k}^2,
\end{equation}
where $\rho_k=a\ip{y^k}{s^k}+b\sigma_k$ and $(a,b)\in\Omega$.

To explain this condition, observe that, by definition of $\rho_k$, the term
$\ip{y^k}{s^k}+\sigma_k$ on the left-hand side of the fundamental
inequality~\eqref{eq:refined_twostep} admits the decomposition
\[
 \ip{y^k}{s^k}+\sigma_k
 =
 (1-a)\ip{y^k}{s^k}+(1-b)\sigma_k+\rho_k.
\]
Since $\ip{y^k}{s^k}\ge0$ and $\sigma_k\ge0$, for $a,b\in[0,1]$ we have
\[
 (1-a)\ip{y^k}{s^k}+(1-b)\sigma_k\ge0.
\]
Dropping this nonnegative quantity from the left-hand side of~\eqref{eq:refined_twostep} gives
\begin{equation*}
\begin{aligned}
 &\norm{x^{k+1}-x^*}^2+\norm{s^{k+1}}^2
 +2t_k(1+\theta_k)\delta_k
 +2\norm{e^k}^2
 +2t_k\theta_k\rho_k\\
 &\qquad\le
 \norm{x^k-x^*}^2+2t_k\theta_k\delta_{k-1}
 +2t_k^2\norm{y^k}^2
 +2t_k(2+\theta_k)\eta_k .
\end{aligned}
\end{equation*}
Move $2t_k\theta_k\rho_k$ from the left-hand side to the right. 
Then, the stepsize condition~\eqref{eq:curvature_condition} turns out to require that the
resulting expression $2t_k^2\norm{y^k}^2+2t_k(2+\theta_k)\eta_k-2t_k\theta_k\rho_k$ be at most
$\norm{s^k}^2$. 
Thus imposing~\eqref{eq:curvature_condition} in the preceding inequality yields
the following estimate corresponding to~\eqref{eq:prior_mm_onestep}:
\[
 \norm{x^{k+1}-x^*}^2+\norm{s^{k+1}}^2
 +2t_k(1+\theta_k)\delta_k
 +2\norm{e^k}^2\le
 \norm{x^k-x^*}^2+\norm{s^k}^2
 +2t_k\theta_k\delta_{k-1}.
\]

When $(a,b)=(1,1)$, the entire term $\ip{y^k}{s^k}+\sigma_k$ can be used to relax the stepsize
condition. 
However, in this setting, the nonnegative term $(1-a)\ip{y^k}{s^k}+(1-b)\sigma_k$ used in the
convergence analysis in Section~\ref{sec:xk_conv} then vanishes and our argument no longer guarantees convergence
of the entire sequence.
For this reason, the parameter region $\Omega$ defined in~\eqref{eq:def_Omega} excludes $(1,1)$.
When $(a,b)=(0,0)$, we have $\rho_k=0$, and, as in the work of Malitsky and
Mishchenko~\cite{MalitskyMishchenko2024}, none of these nonnegative quantities is used to relax
the stepsize condition.

With the iteration data fixed, increasing $a,b$ can enlarge the range of stepsizes
satisfying~\eqref{eq:curvature_condition}.
However, as the numerical experiments in Section~\ref{sec:experiments} illustrate, taking $(a,b)$ closer to $(1,1)$
need not accelerate convergence.
Section~\ref{sec:parameter_variants} examines representative parameter choices and their properties in detail.

We next derive an explicit choice of stepsize satisfying~\eqref{eq:curvature_condition}.
In~\eqref{eq:curvature_condition},
replacing $t_k$ by the unknown $t>0$ and rearranging gives
\begin{equation}\label{eq:stepsize_quadratic}
 2\left(
 \norm{y^k}^2+\frac{\eta_k-\rho_k}{t_{k-1}}
 \right)t^2
 +4\eta_kt-\norm{s^k}^2\le0.
\end{equation}
In what follows, assume $s^k\ne0$.

First suppose that the quadratic coefficient in~\eqref{eq:stepsize_quadratic} is positive.
Since $\eta_k\ge0$, the left-hand side is strictly increasing in $t>0$, and the admissible
positive values of $t$ form an interval whose upper endpoint is the positive root of the quadratic
on the left-hand side.

If the quadratic coefficient is $0$,
then~\eqref{eq:stepsize_quadratic} reduces to
\[
 4\eta_kt\le\norm{s^k}^2.
\]
Thus, if $\eta_k>0$, it suffices that $t\le\norm{s^k}^2/(4\eta_k)$, whereas if $\eta_k=0$, every
$t>0$ is admissible.

If the quadratic coefficient is negative,
the left-hand side of~\eqref{eq:stepsize_quadratic}
need not be monotone in $t>0$, and the positive values of $t$ which satisfy the inequality need
not form a single interval.
In this case we replace the quadratic coefficient by $0$, imposing a stronger condition that
ensures the original inequality.
This gives a uniform description of the admissible positive values of $t$ regardless of the sign
of the quadratic coefficient.

In summary, taking the positive part of the quadratic coefficient and imposing
\begin{equation}\label{eq:stepsize_quadratic_plus}
 2\left[
 \norm{y^k}^2+\frac{\eta_k-\rho_k}{t_{k-1}}
 \right]_+t^2
 +4\eta_kt-\norm{s^k}^2\le0
\end{equation}
ensures that~
\eqref{eq:stepsize_quadratic} holds.
Let $\widehat t_k$ denote the supremum of the positive values of $t$
satisfying~\eqref{eq:stepsize_quadratic_plus}.
If both the quadratic and linear coefficients are $0$, every $t>0$ is admissible, and we set
$\widehat t_k=+\infty$.
Otherwise, solving for the positive root of the polynomial on the left gives
\[
 \widehat t_k=
 \frac{\norm{s^k}^2}{
 2\eta_k+\sqrt{
 4\eta_k^2+
 2\left[
 \norm{y^k}^2+\dfrac{\eta_k-\rho_k}{t_{k-1}}
 \right]_+\norm{s^k}^2
 }}.
\]
This is precisely the stepsize upper bound~\eqref{eq:root} defined in
Section~\ref{sec:stepsize_rule}.

The convergence analysis below also uses the nonnegative terms discarded in deriving the stepsize
condition.
We therefore return to the original fundamental inequality~
\eqref{eq:refined_twostep}.
The stepsize update gives $t_k\le\widehat t_k$, so~\eqref{eq:curvature_condition} holds.
Applying this to~\eqref{eq:refined_twostep}, using the definition of $\rho_k$, and rearranging
yields
\begin{align}
 &\norm{x^{k+1}-x^*}^2+\norm{s^{k+1}}^2
 +2t_k(1+\theta_k)\delta_k
 +2\norm{e^k}^2 +2t_k\theta_k
 \left(
 (1-a)\ip{y^k}{s^k}
 +(1-b)\sigma_k
 \right)\notag\\
 &\qquad\le
 \norm{x^k-x^*}^2+\norm{s^k}^2
 +2t_k\theta_k\delta_{k-1}
 \label{eq:strong_onestep}.
\end{align}

\medskip

Define
\begin{equation}\label{eq:omega}
 \omega_k:=t_k(1+\theta_k)-t_{k+1}\theta_{k+1} \ge 0.
\end{equation}
The nonnegativity of $\omega_k$ follows from $t_{k+1}^2 \le (1+\theta_k)t_k^2$ by the stepsize
update~\eqref{eq:stepsize}.
Using the nonnegativity of $\omega_k$ and~
\eqref{eq:strong_onestep},
we obtain boundedness of the iterates and the following series estimates.

\begin{theorem}[Boundedness of $\{x^k\}$ and $\{s^k\}$]\label{thm:bounded}
Suppose that AdGDC$(a,b)$ does not terminate after finitely many iterations.
For every $x^*\in X^* $, define
\begin{equation}
 D_{x^*}:=
 \norm{x^1-x^*}^2+\norm{s^1}^2
 +2t_1\theta_1\delta_0
 \label{eq:D}.
\end{equation}
Then for every $K\ge1$,
\begin{align}
 &\norm{x^{K+1}-x^*}^2
 +\norm{s^{K+1}}^2
 +2\sum_{k=1}^{K}\omega_k\delta_k
 +2t_{K+1}\theta_{K+1}\delta_K\notag\\
 &\quad
 +2\sum_{k=1}^{K}\norm{e^k}^2
 +2\sum_{k=1}^{K}t_k\theta_k
 \left(
 (1-a)\ip{y^k}{s^k}
 +(1-b)\sigma_k
 \right)
 \le D_{x^*}
 \label{eq:telescoped_rate}.
\end{align}

In particular, $\{x^k\}$ and $\{s^k\}$ are bounded.
\end{theorem}
\begin{proof}
Summing~\eqref{eq:strong_onestep}
over $k=1,\ldots,K$ immediately yields~\eqref{eq:telescoped_rate}.

Since~\eqref{eq:omega} gives $\omega_k\ge0$, \eqref{eq:delta_def} gives $\delta_k\ge0$,
\eqref{eq:ymonotone} gives $\ip{y^k}{s^k}\ge0$, \eqref{eq:prev_gap} gives $\sigma_k\ge0$,
and~\eqref{eq:def_Omega} gives $a,b\in[0,1]$,
every term on the left-hand side of~\eqref{eq:telescoped_rate}
is nonnegative.
Dropping appropriate nonnegative terms gives
\[
 \norm{x^{K+1}-x^*}^2
 +\norm{s^{K+1}}^2
 \le D_{x^*},
\]
so $\{x^k\}$ and $\{s^k\}$ are bounded.
\end{proof}

The estimate in Theorem~\ref{thm:bounded} implies convergence of the series involving the
residuals and correction terms, as stated in the following corollary.

\begin{corollary}\label{cor:summability}
Suppose that AdGDC$(a,b)$ does not terminate after finitely many iterations.
Then
\[
 \sum_{k=1}^{\infty}\norm{e^k}^2<\infty 
\]
and
\begin{equation}\label{eq:summable_curvature}
 \sum_{k=1}^{\infty}t_k\theta_k
 \left(
 (1-a)\ip{y^k}{s^k}
 +(1-b)\sigma_k
 \right)<\infty. 
\end{equation}
\end{corollary}

\begin{proof}
In Theorem~\ref{thm:bounded},
drop the appropriate nonnegative terms from~\eqref{eq:telescoped_rate}
and let $K\to\infty$.
\end{proof}

\section{A lower bound on the stepsizes and the convergence rate}
\label{sec:stepsize_inf_and_conv_rate}
We first show that the stepsizes $t_k$ have a uniform lower bound independent of $k$.
We then use this fact to show that AdGDC$(a,b)$ has a best-iterate convergence rate of $O(1/K)$ in
the number of iterations $K$.
{\bf
Throughout Section~\ref{sec:stepsize_inf_and_conv_rate}, we assume that Algorithm~\ref{alg:adgdc}
does not terminate after finitely many iterations and that the following assumption holds.}

\medskip
\noindent\textbf{Assumption (L).}
For each $F_i\;(i=1,\dots,m)$, $\nabla F_i$ is locally Lipschitz continuous.

\medskip
This condition holds automatically, for example, if each $F_i$ is a convex function of class
$C^2$.

\subsection{A uniform positive lower bound on the stepsizes}
\label{sec:stepsize_inf}

In preparation for bounding $t_k$ from below, the following lemma gives upper bounds on $y^k$ and
$\eta_k$.

\begin{lemma}\label{lem:local_curvature}
For the sequence of iterates, which is bounded by Theorem~\ref{thm:bounded}, there exist finite
$C_y,C_{\eta}\ge0$ such that
\begin{equation}\label{eq:local_curve_bound}
 \norm{y^k}\le C_y\norm{s^k},\qquad
 0\le \eta_k\le C_{\eta}\norm{s^k}^2
 \quad\qquad(k\ge1).
\end{equation}
\end{lemma}
\begin{proof}
Both $x^k$ and $s^k$ are bounded.
Consider a closed ball large enough to contain the iterates.
By continuity, $F$ and $\J$ are bounded on this ball.
By~\eqref{eq:uk}, we have $\widetilde{F}^{k-1}=F(x^{k-1})+\J(x^{k-1})^\top s^k$,
so the bounds established above imply that $\widetilde{F}^{k-1}$ is bounded as well.
Since $g$ is finite-valued and convex, $z^{k-1}\in\partial g(\widetilde{F}^{k-1})$ is uniformly
bounded.
The mapping $\J$ is Lipschitz continuous on a bounded set containing the iterates.
Combining these facts gives
\[\norm{y^k}\le \norm{\J(x^k) -\J(x^{k-1})} \norm{z^{k-1}}\le C_y\norm{s^k}.\]

Moreover, componentwise convexity gives $F(x^k)-\widetilde{F}^{k-1}\ge0$, and a local quadratic
bound on the Taylor remainder gives $\norm{F(x^k)-\widetilde{F}^{k-1}}\le C_F\norm{s^k}^2$.
By Theorem~\ref{thm:subproblem}, we have $z^{k-1}\ge0$.
Together with the componentwise inequality
$F(x^k)-\widetilde{F}^{k-1}\ge0$ established above,
this gives
$\ip{z^{k-1}}{F(x^k)-\widetilde{F}^{k-1}}\ge0$.
Since $g$ is Lipschitz continuous on a bounded set
containing $\widetilde{F}^{k-1}$ and $F(x^k)$,
denoting its Lipschitz constant by $L_g$ yields
\begin{align*}
 0\le \eta_k
 &=g(F(x^k))-g(\widetilde{F}^{k-1})
  -\ip{z^{k-1}}{F(x^k)-\widetilde{F}^{k-1}}\\
 &\le g(F(x^k))-g(\widetilde{F}^{k-1}) \le L_g \norm{F(x^k) - \widetilde{F}^{k-1}}
  \le L_g C_F\norm{s^k}^2 = C_{\eta}\norm{s^k}^2,
\end{align*}
which implies~\eqref{eq:local_curve_bound}.
\end{proof}

Using Lemma~\ref{lem:local_curvature}, we now prove that the stepsizes have a uniform positive
lower bound.

\begin{theorem}[Uniform positive lower bound on the stepsizes]\label{thm:positive_step}
There exists $t_{\inf}>0$ such that
\begin{equation}\label{eq:step_floor}
 t_k\ge t_{\inf}\qquad(k\ge0).
\end{equation}
\end{theorem}
\begin{proof}
For the constants $C_y,C_{\eta}$ in Lemma~\ref{lem:local_curvature}, impose the following
conditions on the lower bound:
\[
 0<t_{\inf}\le t_0,\qquad
 2C_y^2t_{\inf}^2+6C_{\eta} t_{\inf}\le1.
\]
Choose $t_{\inf}$ accordingly.
Assuming $t_{k-1}\ge t_{\inf}$, we obtain
\[
 2t_{\inf}^2\norm{y^k}^2
 +2t_{\inf}\left(2+\frac{t_{\inf}}{t_{k-1}}\right)\eta_k
 \le\left(2C_y^2t_{\inf}^2+6C_{\eta} t_{\inf}\right)\norm{s^k}^2
 \le\norm{s^k}^2.
\]
Here~\eqref{eq:def_rho} gives $\rho_k\ge0$, and~\eqref{eq:eta_k} gives $\eta_k\ge0$.
Therefore,
\[
 \left[
 \norm{y^k}^2+
 \frac{\eta_k-\rho_k}{t_{k-1}}
 \right]_+
 \le
 \norm{y^k}^2+\frac{\eta_k}{t_{k-1}}.
\]
Using the preceding bound yields
\[
 2t_{\inf}^2
 \left[
 \norm{y^k}^2+
 \frac{\eta_k-\rho_k}{t_{k-1}}
 \right]_+
 +4t_{\inf}\eta_k\le
 2t_{\inf}^2\norm{y^k}^2
 +2t_{\inf}
 \left(2+\frac{t_{\inf}}{t_{k-1}}\right)\eta_k
 \le\norm{s^k}^2.
\]
This is the condition obtained by substituting $t=t_{\inf}$
into~\eqref{eq:stepsize_quadratic_plus}.
Hence $\widehat t_k\ge t_{\inf}$ holds.
Moreover,
$\sqrt{1+\theta_{k-1}}\,t_{k-1}\ge t_{k-1}\ge t_{\inf}$
so~\eqref{eq:stepsize} gives $t_k\ge t_{\inf}$.
Induction from the initial value completes the proof.
\end{proof}

\subsection{Best-iterate convergence rate}
\label{sec:conv_rate}

From the estimate in Theorem~\ref{thm:bounded}, we derive a bound on the weighted optimality gaps
that will yield convergence rates for the best iterate and the weighted average.

\begin{lemma}
\label{lem:weighted_gap}
For every $K\ge1$, using $\omega_k$ defined in~\eqref{eq:omega}, define the weights
$\lambda_k^{(K)}$ by
\begin{equation}\label{eq:rate_weights}
 \lambda_k^{(K)}:=
 \begin{cases}
  \omega_k & (1\le k\le K-1)\\
  t_K(1+\theta_K) & (k=K)
 \end{cases}.
\end{equation}
Then $\lambda_k^{(K)}\ge0$, and for every $x^*\in X^* $,
with the constant defined in~\eqref{eq:D},
denoted by $D_{x^*}$, we have
\begin{equation}\label{eq:weighted_gap}
 2\sum_{k=1}^K\lambda_k^{(K)}\delta_k
 \le D_{x^*}.
\end{equation}
Furthermore,
\begin{equation}\label{eq:weight_sum}
 \sum_{k=1}^K\lambda_k^{(K)}
 =\sum_{k=1}^K t_k+t_1\theta_1
 \ge Kt_{\inf}.
\end{equation}
\end{lemma}

\begin{proof}
By $\omega_k\ge0$ and $t_K>0$, we have $\lambda_k^{(K)}\ge0$. 
In Theorem~\ref{thm:bounded}, dropping the appropriate nonnegative terms
from~\eqref{eq:telescoped_rate} gives
\[
 2\left(
 \sum_{k=1}^K\omega_k\delta_k
 +t_{K+1}\theta_{K+1}\delta_K
 \right)\le D_{x^*}.
\]
Since $\omega_K+t_{K+1}\theta_{K+1}=t_K(1+\theta_K)$,
\eqref{eq:weighted_gap} follows.

Moreover, by definition of $\omega_k$,
\[
 \sum_{k=1}^K\lambda_k^{(K)}
 =\sum_{k=1}^{K-1}
 \bigl(t_k(1+\theta_k)-t_{k+1}\theta_{k+1}\bigr)
 +t_K(1+\theta_K)=\sum_{k=1}^K t_k+t_1\theta_1
 \ge Kt_{\inf},
\]
which gives~\eqref{eq:weight_sum}.
\end{proof}

We now use this lemma to establish the best-iterate rate $O(1/K)$.

\begin{theorem}[$O(1/K)$ convergence rate for the best iterate]\label{thm:rate}

Fix $x^*\in X^* $.
With $D_{x^*}$ as in Theorem~\ref{thm:bounded}, for every $K\ge1$ we have
\begin{equation}\label{eq:best_rate}
 \min_{1\le k\le K}\bigl(h(x^k)- h^* \bigr)
 \le\frac{D_{x^*}}{2t_{\inf}K}.
\end{equation}
\end{theorem}
\begin{proof}
By Lemma~\ref{lem:weighted_gap},
\[
 \left(\sum_{k=1}^K\lambda_k^{(K)}\right)
 \min_{1\le k\le K}\delta_k
 \le\sum_{k=1}^K\lambda_k^{(K)}\delta_k
 \le\frac {D_{x^*}}{2}.
\]
Moreover, \eqref{eq:weight_sum} gives
\[
 \min_{1\le k\le K}\delta_k
 \le\frac{D_{x^*}}{2\sum_{k=1}^K\lambda_k^{(K)}}
 \le\frac{D_{x^*}}{2t_{\inf} K},
\]
which proves~\eqref{eq:best_rate}.
\end{proof}

A similar argument gives an $O(1/K)$ rate for the weighted average.

\begin{corollary}
\label{cor:ergodic_rate}
Using the weights defined in~\eqref{eq:rate_weights},
namely $\lambda_k^{(K)}$, define the weighted average of the iterates $\{x^k\}_{k=1}^K$ by
\[
 \bar x^K:=
 \frac{\sum_{k=1}^K\lambda_k^{(K)}x^k}
 {\sum_{k=1}^K\lambda_k^{(K)}}.
\]
Under the assumptions of Theorem~\ref{thm:rate},
\[
 h(\bar x^K)- h^* 
 \le\frac{D_{x^*}}{2t_{\inf}K}.
\]
\end{corollary}

\begin{proof}
Convexity of $h$ and Jensen's inequality give
\[
 h(\bar x^K)- h^* 
 \le
 \frac{\sum_{k=1}^K\lambda_k^{(K)}\delta_k}
 {\sum_{k=1}^K\lambda_k^{(K)}}
 \le
 \frac{D_{x^*}}
 {2\sum_{k=1}^K\lambda_k^{(K)}}
 \le\frac{D_{x^*}}{2t_{\inf}K}.
\]
\end{proof}

\begin{remark}
The $O(1/K)$ bounds in this section apply to the best iterate and the weighted average.
They do not establish $\delta_K=O(1/K)$ for the last iterate.
\end{remark}

\section{Global convergence of the entire sequence}
\label{sec:xk_conv}
Thus far, we have established boundedness of the
iterates, a positive lower bound on the stepsizes,
and a best-iterate convergence rate for AdGDC$(a,b)$.
The rate bounds the difference between the objective value and the optimal value by $O(1/K)$ for
the best iterate $k$ among $1 \le k \le K$.
It remains to show that the objective values $h(x^k)$ along the sequence $\{x^k\}$ converge to the
optimal value $h^*$ and that the entire sequence $\{x^k\}$ converges to a single optimal solution $x^*\in X^*$.
In this section we establish convergence of both the objective values and the iterates themselves.
{\bf
Throughout Section~\ref{sec:xk_conv}, we assume that Algorithm~\ref{alg:adgdc}
does not terminate after finitely many iterations and that Assumption~(L) holds.}

\subsection{Estimates for the residuals and differences}
\label{sec:residual_step}

For the residual $e^k$ defined in~\eqref{eq:ek},
Corollary~\ref{cor:summability} of Section~\ref{sec:meaning_stepsize_ineq} has already established
\[
 \sum_{k=1}^{\infty}\norm{e^k}^2<\infty,
 \qquad
 \sum_{k=1}^{\infty}t_k\theta_k
 \left(
 (1-a)\ip{y^k}{s^k}
 +(1-b)\sigma_k
 \right)<\infty.
\]
We use these estimates to derive the properties of the differences needed to prove convergence of
the entire sequence.

We first bound $\theta_k$ from above using the stepsize update.

\begin{lemma}[Upper bound on $\theta_k$]\label{lem:theta_bound}
Let $\varphi:=(1+\sqrt5)/2$.
Then
\[
 0<\theta_k\le\varphi\qquad(k\ge0).
\]
\end{lemma}
\begin{proof}
We have $\theta_0=1\le\varphi$, and~
\eqref{eq:stepsize} implies
$\theta_k\le\sqrt{1+\theta_{k-1}}$. 
The result follows by induction because $\sqrt{1+\varphi}=\varphi$.
\end{proof}

We next introduce the normalized step direction $v^k$ and relate its successive differences to
$y^k$ and $e^k$.

\begin{lemma}\label{lem:velocity_difference}
For $k\ge0$, let
\begin{equation}\label{eq:vk}
 v^k:=\frac{s^{k+1}}{t_k}.
\end{equation}
Then for every $k\ge1$,
\begin{equation}\label{eq:velocity_difference}
 v^k-v^{k-1}=-y^k+\frac{e^k}{t_k}.
\end{equation}
\end{lemma}
\begin{proof}
Dividing both sides of~\eqref{eq:ek} by $t_k$ gives
\[
 \frac{e^k}{t_k}=\frac{s^{k+1}}{t_k}
 -\frac{s^k}{t_{k-1}}+y^k,
\]
which is~\eqref{eq:velocity_difference}.
\end{proof}

The next lemma bounds $\|y^k\|^2$ in terms of two correction quantities.

\begin{lemma}\label{lem:y_cocoercive}
There exist finite constants $C_{\mathrm{co}},C_\sigma>0$ such that, for every $k\ge1$,
\begin{align}
 \norm{y^k}^2
 &\le C_{\mathrm{co}}\ip{y^k}{s^k},
 \label{eq:cocoercive}\\
 \norm{y^k}^2
 &\le C_\sigma\sigma_k
 \label{eq:sigma_curvature_bound}.
\end{align}
\end{lemma}

\begin{proof}
On a sufficiently large bounded convex set containing the iterates, let $L_i>0$ be a Lipschitz
constant for each $\nabla F_i$.
By uniform boundedness and nonnegativity of $z^{k-1}$, there exists a constant $M_i\ge0$ such that
$0\le z_i^{k-1}\le M_i$.

Set
\begin{align*}
 \Delta_i^k &:=
 \nabla F_i(x^k)-\nabla F_i(x^{k-1}), \\
 B_i^k &:=
 F_i(x^{k-1})-F_i(x^k)
 +\ip{\nabla F_i(x^k)}{s^k}.
\end{align*}
Convexity of $F_i$ gives $B_i^k\ge0$.

By convexity of $F_i$ and Lipschitz continuity of $\nabla F_i$,
\begin{equation}\label{eq:Lsmooth_property}
 F_i(u)-F_i(v)-\ip{\nabla F_i(v)}{u-v}
 \ge \frac{1}{2L_i}
 \norm{\nabla F_i(u)-\nabla F_i(v)}^2.
\end{equation}
Substituting $u=x^{k-1}$ and $v=x^k$ into this inequality yields
\begin{equation}\label{eq:Delta22LB}
  \norm{\Delta_i^k}^2\le 2L_i B_i^k.
\end{equation}
Furthermore, adding the instances of~\eqref{eq:Lsmooth_property} obtained by substituting
$u=x^{k-1}$, $v=x^k$ and $u=x^{k}$, $v=x^{k-1}$, respectively, gives
\begin{equation}\label{eq:Delta2LDs}
   \norm{\Delta_i^k}^2\le L_i\ip{\Delta_i^k}{s^k}.
\end{equation}

On the other hand, the weighted Cauchy--Schwarz inequality and~\eqref{eq:Delta2LDs} give
\begin{align*}
 \norm{y^k}^2
 =\norm{\sum_{i=1}^m z_i^{k-1}\Delta_i^k}^2
 \le
 \left(\sum_{i=1}^m z_i^{k-1}L_i\right)
 \sum_{i=1}^m\frac{z_i^{k-1}}{L_i}
 \norm{\Delta_i^k}^2.
\end{align*}

Hence, by $0\le z_i^{k-1}\le M_i$,
\begin{align*}
 \norm{y^k}^2
 &\le
 \left(\sum_{i=1}^m M_iL_i\right)
 \sum_{i=1}^m z_i^{k-1}
 \ip{\Delta_i^k}{s^k}\\
 &=\left(\sum_{i=1}^m M_iL_i\right)
 \ip{y^k}{s^k}.
\end{align*}

Next, $z^{k-1}\in\partial g(\widetilde F^{k-1})$ implies
\[
 g(F(x^{k-1}))
 \ge g(\widetilde F^{k-1})
 +\ip{z^{k-1}}{F(x^{k-1})-\widetilde F^{k-1}}.
\]
Combining this with the definition~\eqref{eq:eta_k} of $\eta_k$ and the
definition~\eqref{eq:sigma_k} of $\sigma_k$ gives
\begin{align*}
 \sigma_k
 &=g(F(x^{k-1}))-g(\widetilde F^{k-1})-\ip{z^{k-1}}{F(x^k)-\widetilde F^{k-1}}
 +\ip{\J(x^k)z^{k-1}}{s^k}\\
 &\ge
 \ip{z^{k-1}}{
 F(x^{k-1})-F(x^k)+\J(x^k)^\top s^k}\\
 &=\sum_{i=1}^m z_i^{k-1}B_i^k.
\end{align*}

Applying the weighted Cauchy--Schwarz inequality and~\eqref{eq:Delta22LB} once more gives
\begin{align*}
 \norm{y^k}^2
 &\le
 \left(\sum_{i=1}^m z_i^{k-1}L_i\right)
 \sum_{i=1}^m\frac{z_i^{k-1}}{L_i}
 \norm{\Delta_i^k}^2 \\
 &\le
 2\left(\sum_{i=1}^m M_iL_i\right)
 \sum_{i=1}^m z_i^{k-1}B_i^k\\
 &\le
 2\left(\sum_{i=1}^m M_iL_i\right)\sigma_k.
\end{align*}

Therefore, setting
\[
 C_{\mathrm{co}}:=
 \max\left\{1,\;\sum_{i=1}^mM_iL_i\right\},
 \qquad
 C_\sigma:=
 \max\left\{1,\;2\sum_{i=1}^mM_iL_i\right\}
\]
proves the claim.
\end{proof}

We use this estimate to establish convergence of a weighted series involving $\|y^k\|^2$.
\begin{lemma}
\label{lem:y_weighted_summability}
The following series converges:
\begin{equation}\label{eq:y_squared_weighted}
 \sum_{k=1}^{\infty}
 t_k\theta_k\norm{y^k}^2<\infty.
\end{equation}
\end{lemma}

\begin{proof}
By Lemma~\ref{lem:y_cocoercive}, for every $k\ge1$,
\[
 (1-a)\ip{y^k}{s^k}+(1-b)\sigma_k
 \ge
 \left(
 \frac{1-a}{C_{\mathrm{co}}}
 +\frac{1-b}{C_\sigma}
 \right)\norm{y^k}^2.
\]

Since $(a,b)\in\Omega$ implies both $a,b\in[0,1]$ and $a+b<2$,
\[
 C_{a,b}:=
 \frac{1-a}{C_{\mathrm{co}}}
 +\frac{1-b}{C_\sigma}>0.
\]
Hence Corollary~\ref{cor:summability},
specifically~\eqref{eq:summable_curvature}, gives
\[
\begin{aligned}
 C_{a,b}
 \sum_{k=1}^{\infty}t_k\theta_k\norm{y^k}^2
 &\le
 \sum_{k=1}^{\infty}t_k\theta_k
 \left(
 (1-a)\ip{y^k}{s^k}
 +(1-b)\sigma_k
 \right)\\
 &<\infty .
\end{aligned}
\]
The claim follows from $C_{a,b}>0$.
\end{proof}

\subsection{Convergence of the objective values and of the entire sequence}

The following lemma shows that the objective values along the iterates converge to the optimal
value.

\begin{lemma}[Convergence of the objective values]\label{lem:objective_convergence}
For the sequence of iterates $\{x^k\}$,
\[
 \delta_k=h(x^k)- h^* \longrightarrow0
 \qquad(k\to\infty).
\]
\end{lemma}
The proof is given in Appendix~\ref{app:convergence_only}.

\medskip
Finally, we show that the sequence $\{x^k\}$ generated by AdGDC$(a,b)$ has a limit in the solution
set $X^*$.

\begin{theorem}[Global convergence of the entire sequence]\label{thm:global}
There exists $x^* \in X^* $ such that
\[
 x^k\longrightarrow x^* \qquad(k\to\infty).
\]
In particular, $s^k\to0$.
\end{theorem}

The proof is given in Appendix~\ref{app:convergence_only}.

\section{Representative parameter choices and their properties}
\label{sec:parameter_variants}

We examine representative parameter choices for AdGDC$(a,b)$, comparing the correction terms they
use and studying how the parameters affect the stepsize update.
The convergence results established above hold for every $(a,b)\in\Omega$, so we focus here on the
correction terms and the properties of the parameters.

\subsection{Representative choices and comparison of the correction terms}

The stepsize condition for AdGDC$(a,b)$ uses $y^k$ and $\sigma_k$, defined in~\eqref{eq:yk}
and~\eqref{eq:sigma_k}, to form the correction term
\[
 \rho_k=a\ip{y^k}{s^k}+b\sigma_k.
\]
Here $y^k$ is a vector formed from the difference of the Jacobians at successive iterates, and
$\sigma_k$ is a nonnegative quantity defined in terms of the error in a first-order approximation
of the objective.
Representative parameter choices include
\begin{align*}
 (a,b)=(0,0)&:\quad \rho_k=0,\\
 (a,b)=(1,0)&:\quad \rho_k=\ip{y^k}{s^k},\\
 (a,b)=(0,1)&:\quad \rho_k=\sigma_k.
\end{align*}

The baseline choice $(a,b)=(0,0)$ uses no correction term and determines the stepsize from
$\norm{y^k}^2$ and the first-order approximation error $\eta_k$ of the outer function $g$.
The choice $(a,b)=(1,0)$ uses the inner product $\ip{y^k}{s^k}$ between the change in the
Jacobian, with the prescribed weights, and the displacement between iterates.
When $g$ is affine, this quantity corresponds to the inner product involving the gradient
difference in smooth convex optimization.
The choice $(a,b)=(0,1)$ uses $\sigma_k$, which incorporates a first-order approximation error
arising from the composite structure. 
Moreover, taking $a,b>0$ allows both correction terms to be used together.

Both $\ip{y^k}{s^k}$ and $\sigma_k$ are nonnegative quantities arising from convexity, but
neither dominates the other in general. 
The following proposition compares them by giving an explicit expression for their difference.

\begin{proposition}
\label{prop:correction_comparison}
For every $k\ge1$,
\begin{equation}\label{eq:correction_comparison}
\begin{aligned}
 \sigma_k-\ip{y^k}{s^k}=
 \Bigl[
 g(F(x^{k-1}))-g(\widetilde F^{k-1})&
 -\ip{z^{k-1}}{
 F(x^{k-1})-\widetilde F^{k-1}}
 \Bigr]\\
 &\qquad\quad
 -\ip{z^{k-1}}{
 F(x^k)-\widetilde F^{k-1}}.
\end{aligned}
\end{equation}
Both the quantity in square brackets and the last inner product on the right-hand side are
nonnegative.

In particular, if $g$ is affine, then
\[
 0\le\sigma_k\le\ip{y^k}{s^k},
\]
whereas if $F$ is affine, then
\[
 0=\ip{y^k}{s^k}\le\sigma_k.
\]
\end{proposition}

\begin{proof}
The definitions of $y^k$ and $\sigma_k$ give
\[
 \sigma_k-\ip{y^k}{s^k}
 =h(x^{k-1})-h(x^k)+\eta_k+\ip{\J(x^{k-1})z^{k-1}}{s^k}.
\]
Substituting the definition of $\eta_k$ and the identity
\[
 \widetilde F^{k-1}
 =F(x^{k-1})+\J(x^{k-1})^\top s^k
\]
into this expression yields~
\eqref{eq:correction_comparison}.

Moreover, the inclusion
$z^{k-1}\in\partial g(\widetilde F^{k-1})$
implies
\[
 g(F(x^{k-1}))-g(\widetilde F^{k-1})
 -\ip{z^{k-1}}{
 F(x^{k-1})-\widetilde F^{k-1}}
 \ge0.
\]
By convexity of each $F_i$, we also have the componentwise inequality
\[
 F(x^k)\ge
 F(x^{k-1})+\J(x^{k-1})^\top s^k
 =\widetilde F^{k-1}.
\]
Since $z^{k-1}\ge0$ by Theorem~\ref{thm:subproblem}, we obtain
\[
 \ip{z^{k-1}}{
 F(x^k)-\widetilde F^{k-1}}
 \ge0.
\]

If $g$ is affine, then in~
\eqref{eq:correction_comparison}
the bracketed quantity is $0$, giving $\sigma_k\le\ip{y^k}{s^k}$.
If $F$ is affine, then $y^k=0$ and $F(x^k)=\widetilde F^{k-1}$, giving
$\ip{y^k}{s^k}=0\le\sigma_k$.
\end{proof}

By Proposition~\ref{prop:correction_comparison}, the difference between the correction terms
equals the first-order approximation error of the outer function $g$ minus the first-order
approximation error of the inner mapping $F$ weighted by $z^{k-1}$.
In particular, when $g$ is affine, $\ip{y^k}{s^k}$ is at least $\sigma_k$, whereas when $F$ is
affine the reverse inequality holds.
Thus for general nonlinear composite problems there is no uniform ordering of the two correction
terms, and each captures different information.

\subsection{Properties of the parameter choices}

We next examine how $(a,b)$ affects stepsize selection.
At iteration $k$, hold
$s^k,y^k,\eta_k,\sigma_k,t_{k-1},\theta_{k-1}$
fixed and suppose that $s^k\ne0$.
Denote the correction term for $(a,b)$ by $\rho_k(a,b)$.
The stepsize upper bound is
\begin{equation}\label{eq:parameter_hat_t}
 \widehat t_k(a,b)=
 \frac{\norm{s^k}^2}{
 2\eta_k+
 \sqrt{
 4\eta_k^2+
 2\left[
 \norm{y^k}^2+
 \dfrac{\eta_k-\rho_k(a,b)}{t_{k-1}}
 \right]_+
 \norm{s^k}^2
 }
 },
\end{equation}
with the convention $\widehat t_k(a,b)=+\infty$ when the denominator is $0$.

The next proposition shows that increasing the parameters cannot decrease the stepsize upper
bound.

\begin{proposition}
\label{prop:parameter_monotonicity}
Suppose that $(a,b),(a',b')\in\Omega$ satisfy
\[
 a\le a',\qquad b\le b'.
\]
Then
\begin{equation}\label{eq:parameter_monotonicity}
 \widehat t_k(a,b)
 \le\widehat t_k(a',b').
\end{equation}
\end{proposition}

\begin{proof}
By $\ip{y^k}{s^k}\ge0$ and $\sigma_k\ge0$,
\[
 \rho_k(a',b')-\rho_k(a,b)=
 (a'-a)\ip{y^k}{s^k}
 +(b'-b)\sigma_k
 \ge0.
\]
Hence
\[
 \left[
 \norm{y^k}^2+
 \frac{\eta_k-\rho_k(a',b')}{t_{k-1}}
 \right]_+\le
 \left[
 \norm{y^k}^2+
 \frac{\eta_k-\rho_k(a,b)}{t_{k-1}}
 \right]_+.
\]
By~\eqref{eq:parameter_hat_t},
$\widehat t_k$ is nonincreasing in this positive part, which gives~
\eqref{eq:parameter_monotonicity}.
\end{proof}

Even a strict increase in $\widehat t_k$ need not increase the actual stepsize $t_k$.
Indeed, since
\[
 t_k=
 \min\left\{
 \sqrt{1+\theta_{k-1}}\,t_{k-1},
 \widehat t_k
 \right\},
\]
an increase in $\widehat t_k$ leaves $t_k$ unchanged whenever the growth restriction is active.
Moreover, Proposition~\ref{prop:parameter_monotonicity} compares the bounds for fixed iteration
data.
It does not imply an ordering of the iterates or convergence speeds obtained with different parameter settings.

Finally, consider the parameter region
\[
 \Omega=\{(a,b)\in[0,1]^2:a+b<2\}.
\]
We defined this region as above because it ensures the estimates required in the convergence
analysis.

Indeed, the estimate~\eqref{eq:strong_onestep} contains
\[
 (1-a)\ip{y^k}{s^k}+(1-b)\sigma_k
\]
as a nonnegative term.
Moreover, by Lemma~\ref{lem:y_cocoercive},
specifically~\eqref{eq:cocoercive} and~
\eqref{eq:sigma_curvature_bound},
local Lipschitz continuity of each $\nabla F_i$ and boundedness of the iterates imply that there
exist constants $C_{\mathrm{co}},C_\sigma>0$ such that
\[
 \norm{y^k}^2
 \le C_{\mathrm{co}}\ip{y^k}{s^k},
 \qquad
 \norm{y^k}^2
 \le C_\sigma\sigma_k.
\]
Hence
\[
(1-a)\ip{y^k}{s^k}+(1-b)\sigma_k\ge
 \left(
 \frac{1-a}{C_{\mathrm{co}}}
 +\frac{1-b}{C_\sigma}
 \right)\norm{y^k}^2.
\]

If $(a,b)\in\Omega$, then
\[
 \frac{1-a}{C_{\mathrm{co}}}
 +\frac{1-b}{C_\sigma}>0,
\]
and this positive coefficient yields convergence of the weighted series involving $\{y^k\}$. 
The argument requires a bound on $\norm{y^k}^2$ with a strictly positive coefficient in addition
to nonnegativity of the term.

At $(a,b)=(1,1)$ the coefficient is $0$, so our convergence proof does not apply directly. 
This does not imply that the algorithm fails to converge for this choice. 
Indeed, for every $\epsilon\in(0,1]$, our convergence guarantee holds at
$(a,b)=(1-\epsilon,1-\epsilon)$. 
Convergence for $(a,b)=(1,1)$ is not established by the present analysis and is left for future
investigation.

\section{Numerical experiments}
\label{sec:experiments}

We evaluate AdGDC$(a,b)$ on nonsmooth convex minimization problems with a maximum-type objective.
First, we compare AdGDC$(a,b)$ for several choices of $(a,b)$ with the methods of Drusvyatskiy and
Paquette~\cite{DrusvyatskiyPaquette2019}.
We then vary the correction coefficients $(a,b)$ in a grid search to examine how the parameters
$(a,b)$ affect convergence.

\subsection{Problem formulation and experimental setup}
\label{subsec:experiments_setting}

We construct the following convex optimization problem based on the log-cosh loss:
\begin{equation}\label{eq:exp_penalty_objective}
 \min_{x\in\R^n}\quad h(x)
 :=\max_{1\le i\le m}F_i(x),
\end{equation}
where
\begin{equation}
 F_i(x):=
 \frac{\lambda}{2}\norm{x}^2
 +\frac{1}{p}\sum_{j=1}^{p}
 \log\cosh\bigl(a_{ij}^{\top}x-y_{ij}\bigr)
 \qquad(i=1,\ldots,m)\notag,
\end{equation}
and the regularization coefficient is $\lambda=0.5$.
We set $n=50$, $m=150$, and $p=20$.
Taking $F(x)=(F_1(x),\ldots,F_m(x))$ and $g(u)=\max_{1\le i\le m}u_i$ puts this problem in the
form~\eqref{eq:problem_intro}.

We draw the vectors $a_{ij}\in\R^n$ independently from the multivariate normal distribution
$\mathcal N(0,I_n/2)$ and randomly generate a vector $x^\dagger\in\R^n$ satisfying
$\norm{x^\dagger}=2$. 
The index $i$ identifies the data group and $j$ the observation within that group. 
Each group comprises $p$ observations.~\footnote{
The underlying regression problem is to estimate an unknown coefficient vector $x^\dagger$ from
observed data $(a_{ij},y_{ij})$. 
Here $a_{ij}$ is the predictor vector, $y_{ij}$ is the response, and the decision variable $x$
represents the coefficient vector to be estimated. 
In contrast to ordinary empirical loss minimization, minimizing the maximum groupwise average loss
emphasizes the group with the largest loss.
}
We generate the observations $y_{ij}$ according to
\[
 y_{ij}:=
 a_{ij}^{\top}x^\dagger+
 \exp(\nu\zeta_i)\varepsilon_{ij},
\]
where $\zeta_i$ and $\varepsilon_{ij}$ are mutually independent standard normal random variables.
Observations within a group share the noise standard deviation $\exp(\nu\zeta_i)$.
Larger values of $\nu$ lead to greater variation in noise standard deviations across groups.
We consider two settings, $\nu=0.3$ and $\nu=1.5$.

For each $\nu\in\{0.3, 1.5\}$, we generate 10 instances as above, using common values of $a_{ij}$,
$x^\dagger$, $\zeta_i$, and $\varepsilon_{ij}$. 
For each instance we generate ten random initial points satisfying $\norm{x^0}=1$ and run all
methods from the same initial points. 
Thus for each $\nu$ we compare 100 runs, comprising 10 instances with 10 starting points each.

We test the proposed method with $(a,b)=(0,0)$, $(1,0)$, $(0,1)$, and $(0.9,0.9)$. 
The first three are the representative choices discussed in Section~\ref{sec:parameter_variants},
and $(0.9,0.9)$ assigns large weights to both correction terms. 
We also include the prox-linear method~\cite[Algorithm~1]{DrusvyatskiyPaquette2019} and
accelerated prox-linear method~\cite[Algorithm~7]{DrusvyatskiyPaquette2019} proposed by
Drusvyatskiy and Paquette. 
Their update rules are given in Appendix~\ref{app:prox_linear}.

As an upper bound on the Lipschitz constant of the gradient of each $F_i$, we use
\begin{equation}\label{eq:exp_L}
 L_i:=\lambda+\frac{\norm{A_i}^2}{p},
 \qquad
 L:=\max_{1\le i\le m}L_i,
\end{equation}
where $A_i\in\R^{p\times n}$ is the matrix whose $j$th row is $a_{ij}^{\top}$.
Indeed, $(\log\cosh r)''\le1$ gives
\[
 \nabla^2F_i(x)
 \preceq
 \lambda I+\frac{1}{p}A_i^\top A_i.
\]
For all variants of AdGDC$(a,b)$, we initialize with $t_0=1/L$ and $\theta_0=1$. 
We use this initialization to give the methods comparable starting conditions. 
Note that the AdGDC$(a,b)$ update itself requires no prior estimate of $L$. 
For the prox-linear method, we fix the stepsize at $t=1/L$. 
For the accelerated prox-linear method, we set the coefficient of the proximal regularization term
to $\widetilde\mu=1.01L$. 
Neither comparison method uses a line search.

At each iteration we dualize the subproblem to obtain a convex quadratic program over the simplex
and solve it with OSQP.
We set both the absolute and relative tolerances to $10^{-9}$ and enable polishing.

We use the stopping criterion
\begin{equation}\label{eq:exp_stop}
 h(x^k)-h^*\le10^{-8},
\end{equation}
and cap the number of iterations at 500.
We compute the reference optimal value $h^*$ separately for the experiments.
It is used only in the stopping test and performance evaluation.
The iterates and stepsizes of all methods are independent of this reference value.

Our primary performance measure is the number of iterations until the stopping criterion is first
satisfied.
Each method solves one quadratic program per iteration, so the number of subproblem solves to
convergence equals the iteration count.

\subsection{Results}
\label{subsec:experiments_results}

Table~\ref{tab:exp_method_summary} reports the mean and median iteration counts for AdGDC$(a,b)$
with the representative parameter choices and for the two comparison methods, for each value of
$\nu$.

\begin{table}[H]
 \centering
\caption{Comparison of iteration counts.
Each value is based on 100 runs.}
 \label{tab:exp_method_summary}
 \begin{tabular}{lrrrr}
  \toprule
  & \multicolumn{2}{c}{$\nu=0.3$}
  & \multicolumn{2}{c}{$\nu=1.5$}\\
  \cmidrule(lr){2-3}\cmidrule(lr){4-5}
Method & Mean & Median & Mean & Median\\
  \midrule
  AdGDC$(0,0)$       & 13.98 & 14.0 & 13.46 & 13.0\\
  AdGDC$(1,0)$       & 11.48 & 11.0 & 12.61 & 12.0\\
  AdGDC$(0,1)$       & 11.98 & 12.0 & 10.52 & 11.0\\
  AdGDC$(0.9,0.9)$   & 13.43 & 13.0 & 14.83 & 14.5\\
  Prox-linear~\cite{DrusvyatskiyPaquette2019}        & 64.34 & 64.0 & 77.05 & 77.0\\
  Accelerated prox-linear~\cite{DrusvyatskiyPaquette2019}
                       & 58.20 & 59.0 & 55.81 & 56.0\\
  \bottomrule
 \end{tabular}
\end{table}

For each $\nu$, all methods met the stopping criterion~\eqref{eq:exp_stop} in all 100 runs.
At $\nu=0.3$, AdGDC$(1,0)$ had the lowest mean iteration count (11.48), followed by AdGDC$(0,1)$
(11.98).
At $\nu=1.5$, AdGDC$(0,1)$ had the lowest mean (10.52), while AdGDC$(1,0)$ required 12.61
iterations on average.
Thus the relative performance of AdGDC$(1,0)$ and AdGDC$(0,1)$, which use the two correction terms
separately, reversed as $\nu$ changed.

\begin{figure}[H]
 \centering
 \begin{minipage}[c]{.33\linewidth}
  \centering
  \includegraphics[width=1.0\linewidth]
   {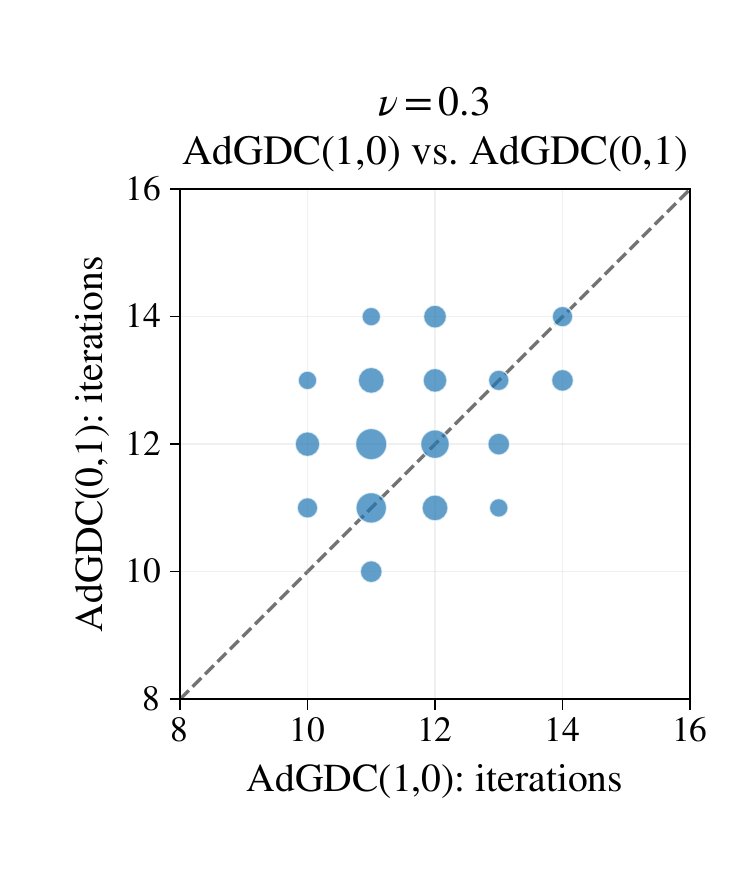}
 \end{minipage}\hfill
 \begin{minipage}[c]{.33\linewidth}
  \centering
  \includegraphics[width=1.0\linewidth]
   {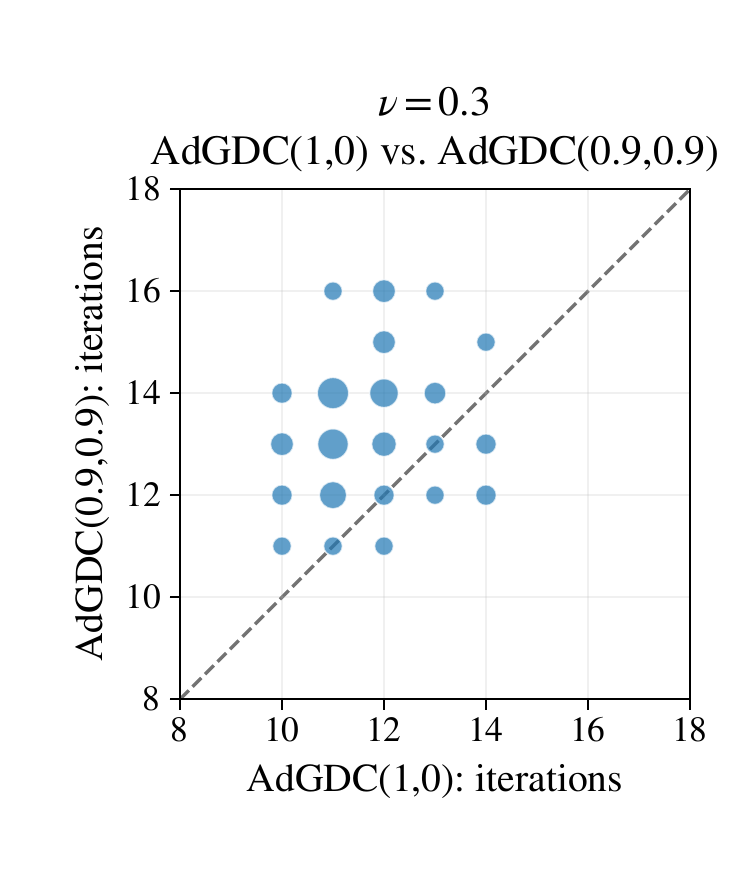}
 \end{minipage}\hfill
 \begin{minipage}[c]{.33\linewidth}
  \centering
  \includegraphics[width=1.0\linewidth]
   {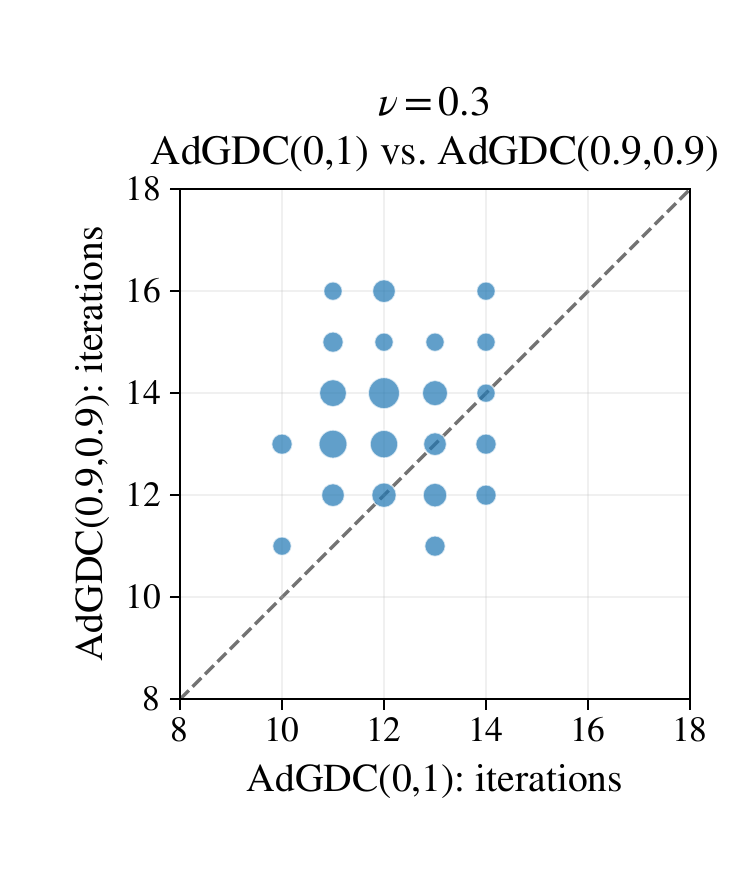}
 \end{minipage}
\caption{Pairwise comparison of iteration counts for the $(a,b)$ settings ($\nu=0.3$).}
 \label{fig:exp_scatter_low}
\end{figure}

\begin{figure}[H]
 \centering
 \begin{minipage}[c]{.33\linewidth}
  \centering
  \includegraphics[width=1.0\linewidth]
   {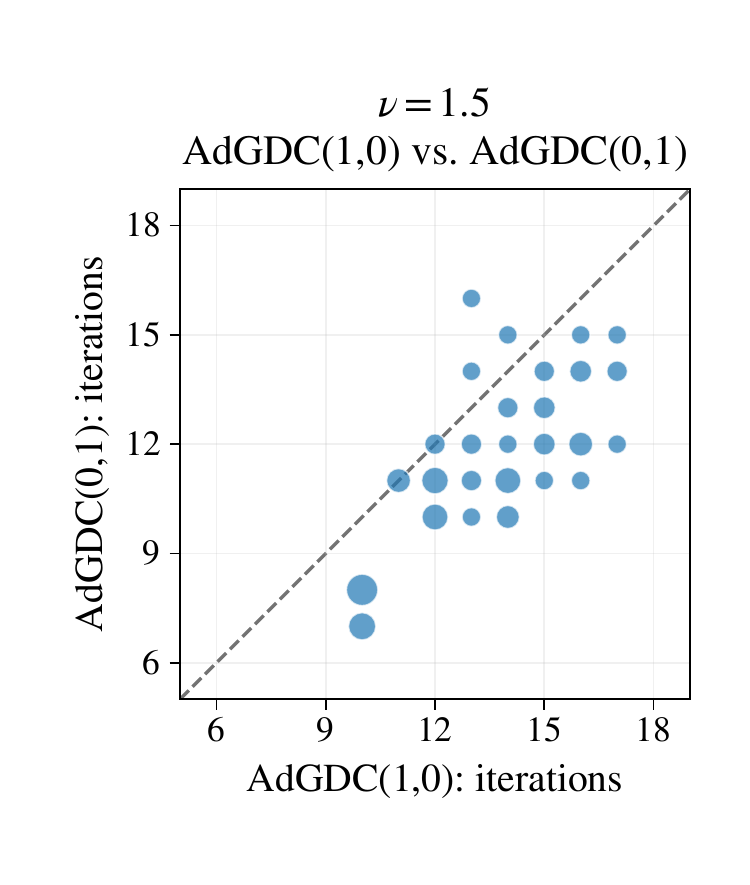}
 \end{minipage}\hfill
 \begin{minipage}[c]{.33\linewidth}
  \centering
  \includegraphics[width=1.0\linewidth]
   {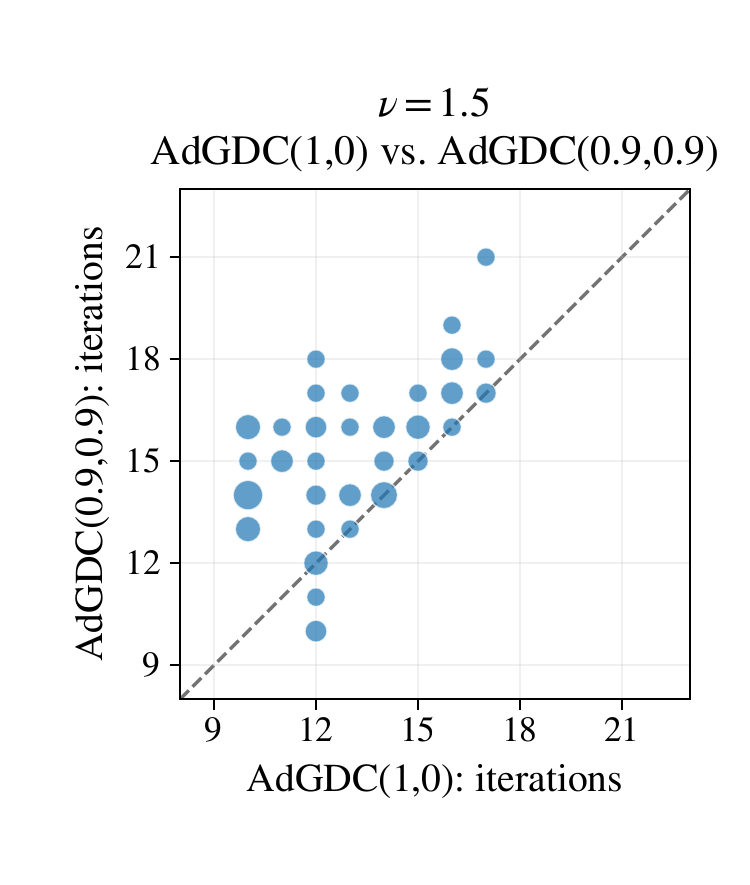}
 \end{minipage}\hfill
 \begin{minipage}[c]{.33\linewidth}
  \centering
  \includegraphics[width=1.0\linewidth]
   {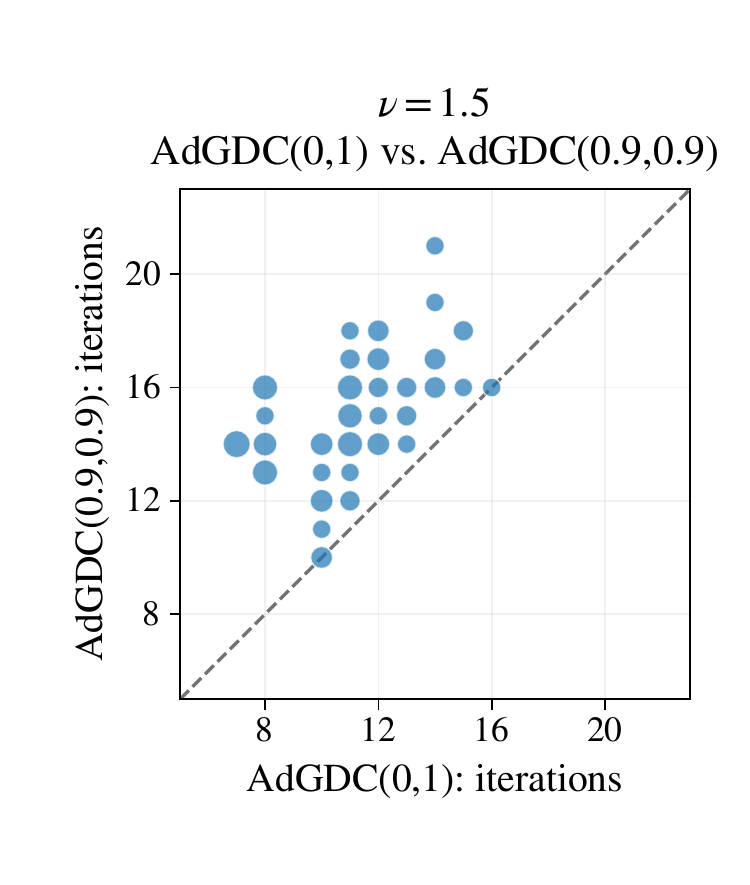}
 \end{minipage}
\caption{Pairwise comparison of iteration counts for the $(a,b)$ settings ($\nu=1.5$).}
 \label{fig:exp_scatter_high}
\end{figure}

To examine this difference more closely, Figures~\ref{fig:exp_scatter_low}
and~\ref{fig:exp_scatter_high} compare iteration counts on matched instances and starting points.
Each figure compares AdGDC$(1,0)$ with AdGDC$(0,1)$, AdGDC$(1,0)$ with AdGDC$(0.9,0.9)$, and
AdGDC$(0,1)$ with AdGDC$(0.9,0.9)$.
Points on the diagonal indicate equal counts.

Comparing AdGDC$(1,0)$ with AdGDC$(0,1)$ at $\nu=0.3$, the former required fewer iterations in 47
runs, the latter in 18, and the two tied in the remaining 35. 
At $\nu=1.5$, AdGDC$(1,0)$ required fewer iterations in 3 runs, AdGDC$(0,1)$ in 90, and they tied
in the remaining 7. 
Thus the reversal in the iteration counts in Table~\ref{tab:exp_method_summary} reflects a pattern
across many instances and starting points, rather than a large difference in only a few runs.

\begin{figure}[H]
 \centering
 \includegraphics[width=.49\linewidth]
  {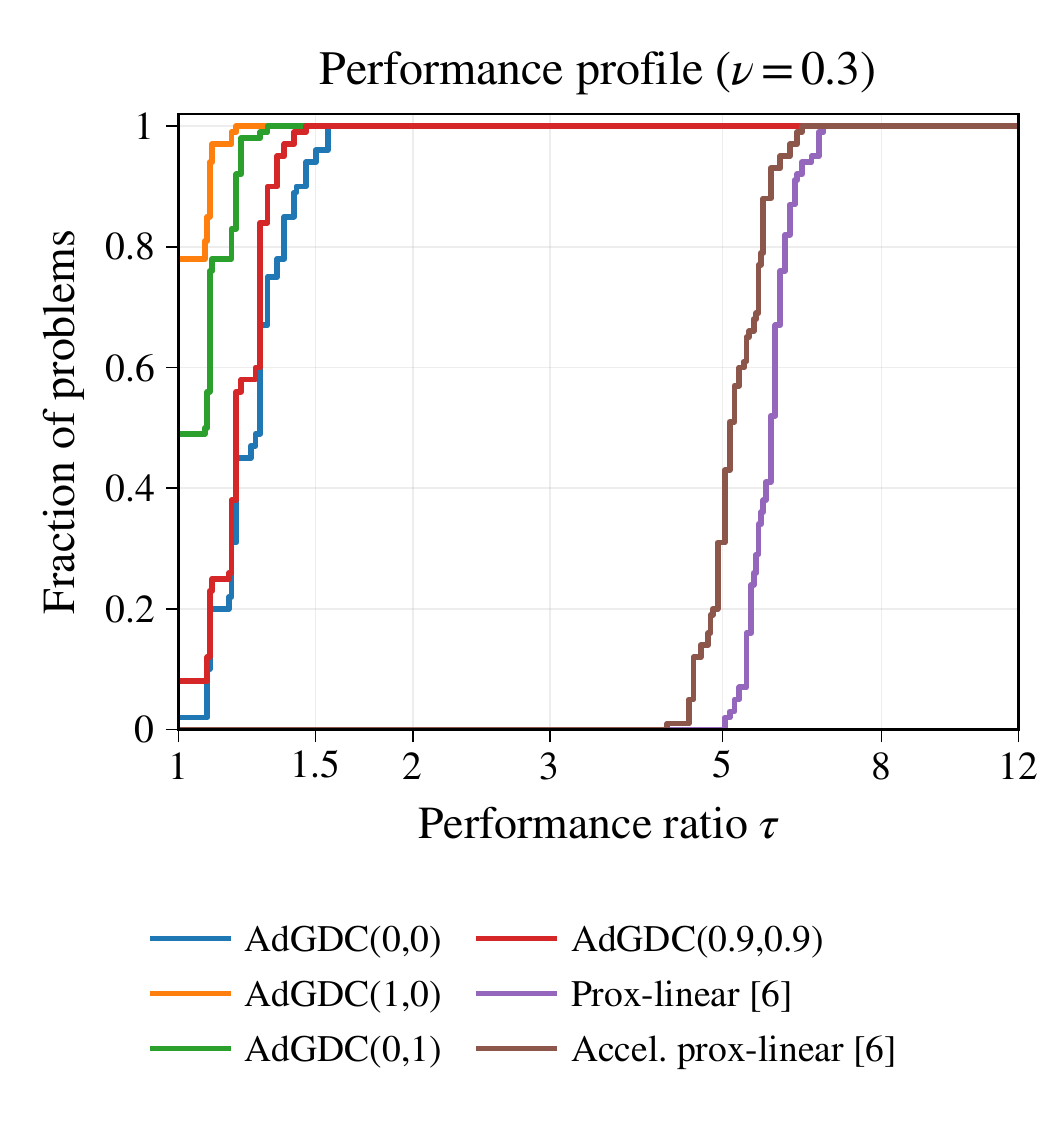}
 \hfill
 \includegraphics[width=.49\linewidth]
  {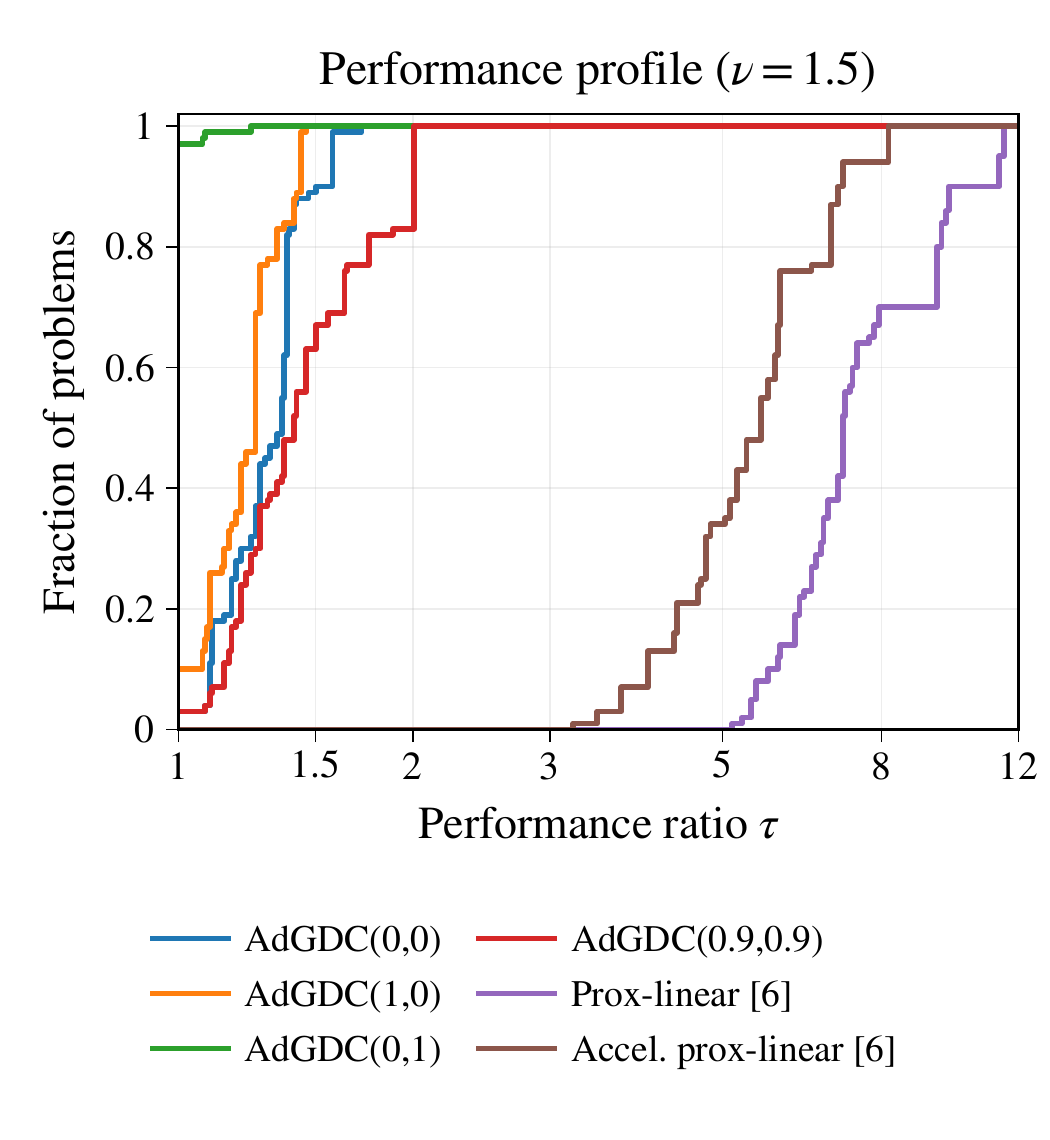}
\caption{Performance profiles based on iteration counts.
Left: $\nu=0.3$.
Right: $\nu=1.5$.}
 \label{fig:exp_profiles}
\end{figure}

We also compare the performance profiles of Dolan and Mor\'e~\cite{DolanMore2002}, dividing each
method's iteration count in each run by the minimum count among the six methods in that run.
Figure~\ref{fig:exp_profiles} shows the cumulative fraction of runs for which this ratio is at
most $\tau$.

Figure~\ref{fig:exp_profiles} likewise shows a change in the relative positions of AdGDC$(1,0)$
and AdGDC$(0,1)$ between $\nu=0.3$ and $\nu=1.5$.
The prox-linear and accelerated prox-linear methods tended to require more iterations than
AdGDC$(a,b)$ in both settings.

\begin{figure}[H]
 \centering
 \includegraphics[width=.49\linewidth]
  {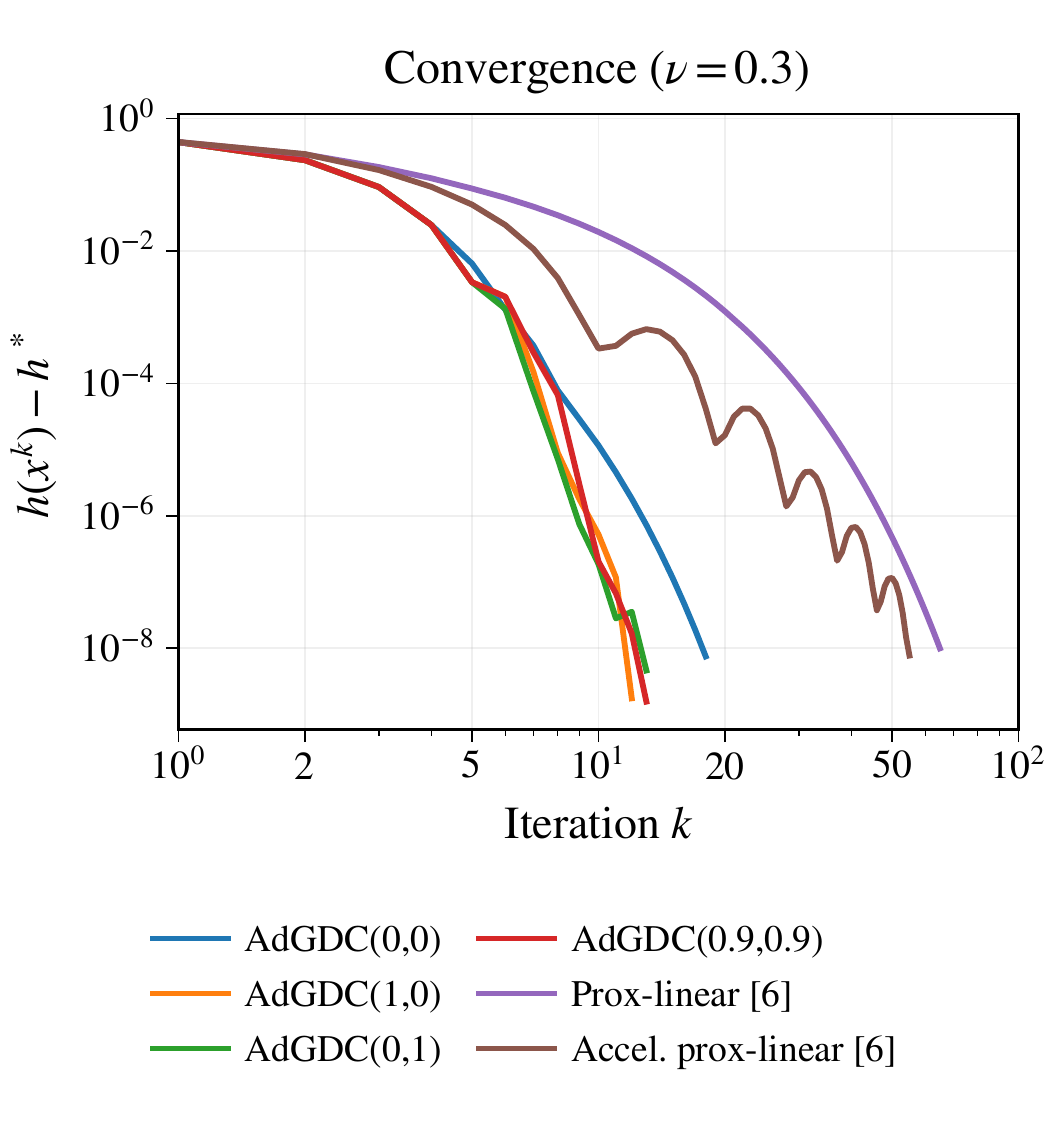}
 \hfill
 \includegraphics[width=.49\linewidth]
  {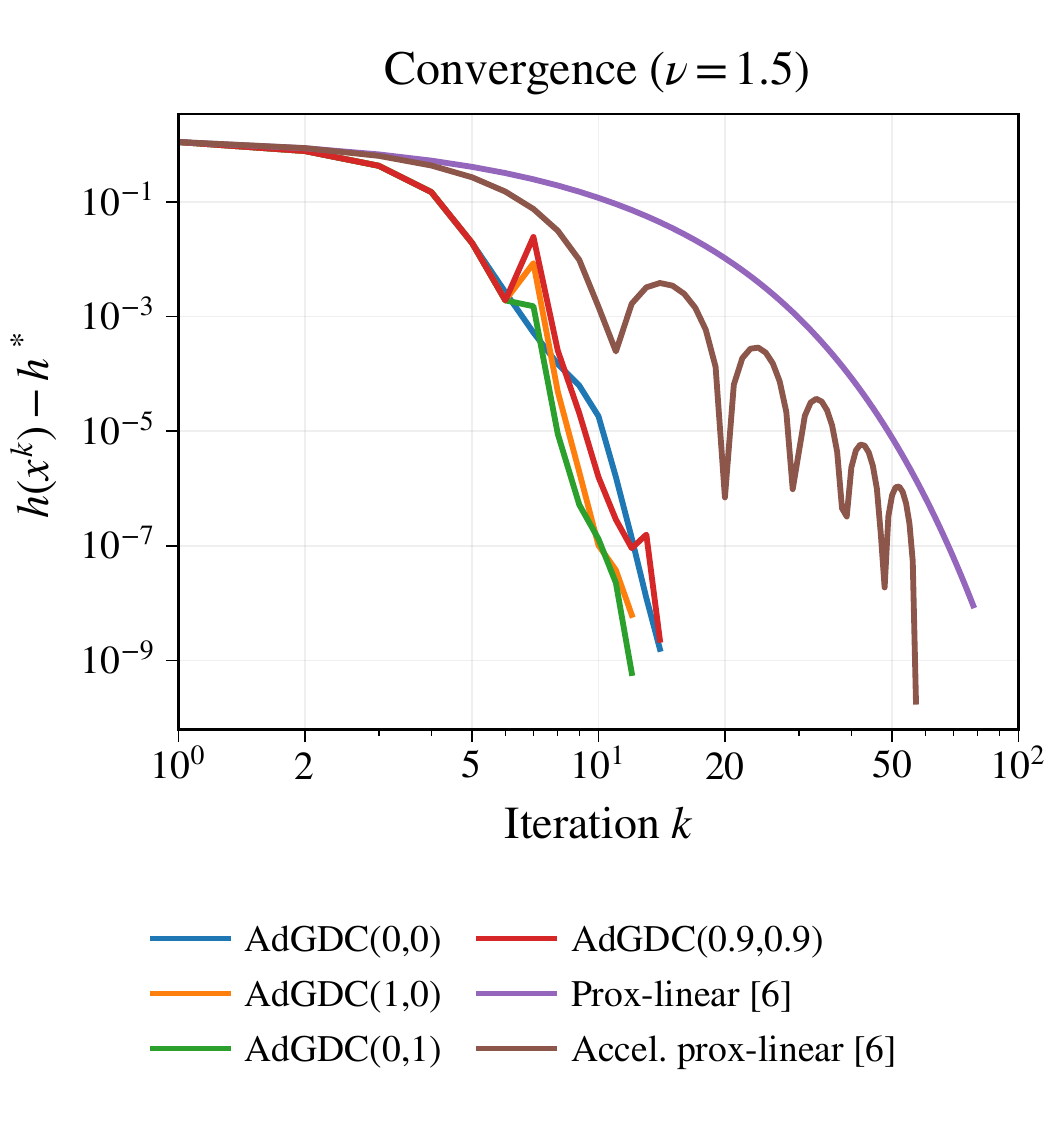}
\caption{Evolution of the optimality gap
$h(x^k)-h^*$ in a representative run.
Left: $\nu=0.3$.
Right: $\nu=1.5$.}
 \label{fig:exp_convergence}
\end{figure}

\begin{figure}[H]
 \centering
 \includegraphics[width=.49\linewidth]
  {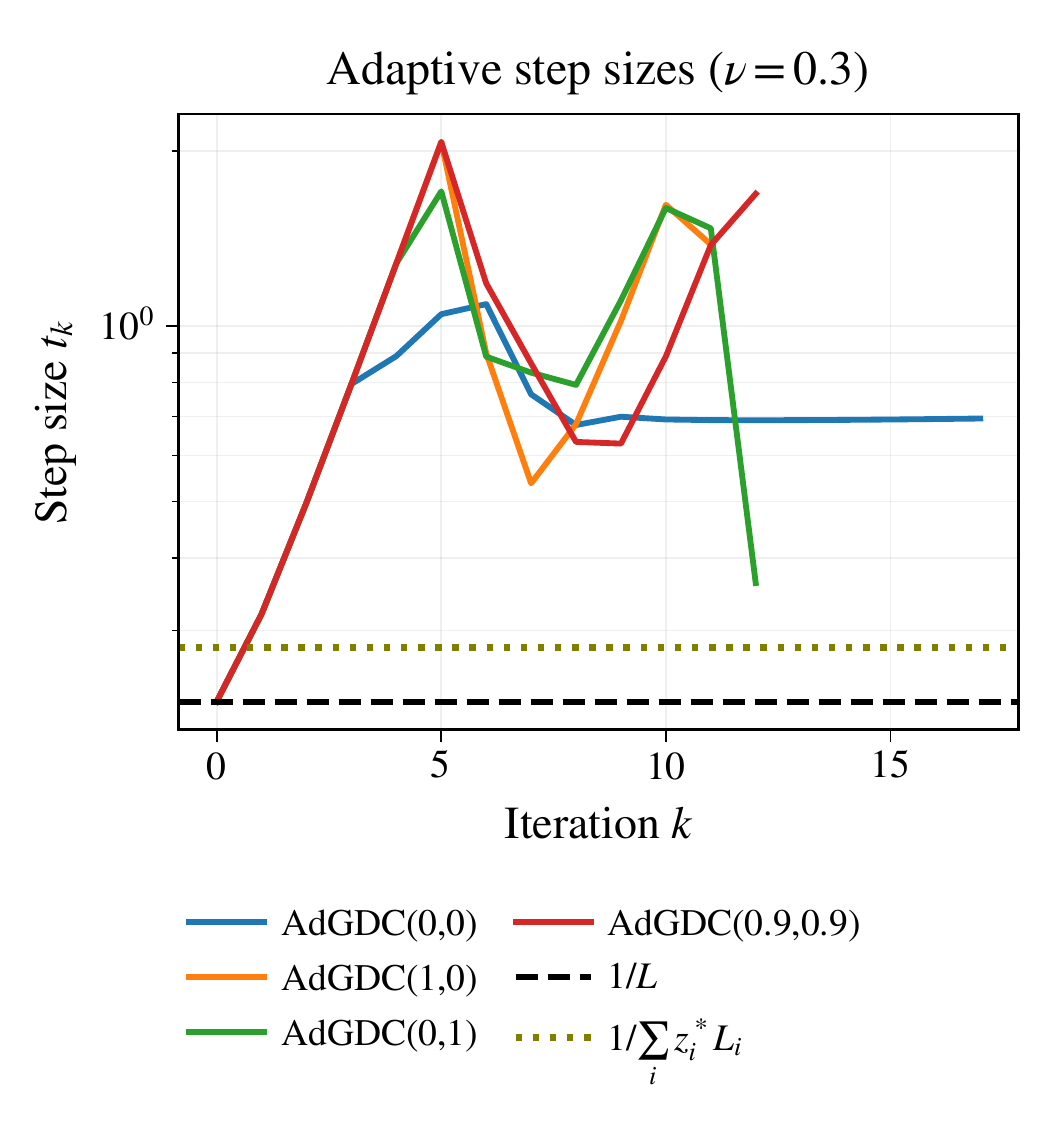}
 \hfill
 \includegraphics[width=.49\linewidth]
  {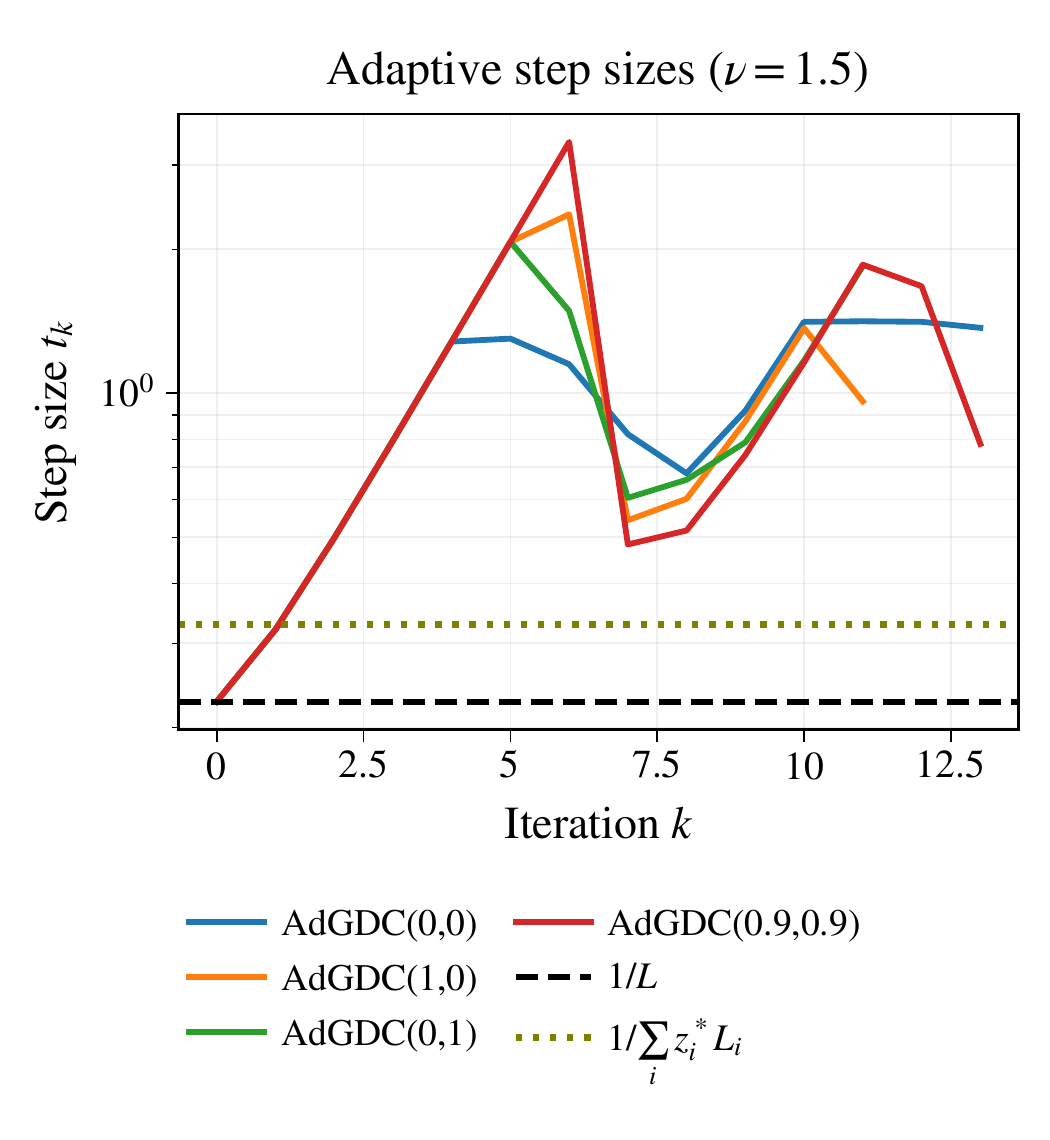}
\caption{Evolution of the AdGDC$(a,b)$
stepsize $t_k$ in a representative run.
Left: $\nu=0.3$.
Right: $\nu=1.5$.
The reference values $1/L$ and $1/\sum_{i=1}^m z_i^*L_i$ are also shown.
Here $z_i^*$ is the dual variable corresponding to the reference optimal solution.}
 \label{fig:exp_stepsizes}
\end{figure}

For one representative run, Figure~\ref{fig:exp_convergence} shows the evolution of the optimality
gap $h(x^k)-h^*$.
Figure~\ref{fig:exp_stepsizes} shows the stepsizes $t_k$ of AdGDC$(a,b)$ on the same instance.

Figure~\ref{fig:exp_stepsizes} includes two benchmarks: $1/L$, based on the global Lipschitz
constant, and $1/\sum_{i=1}^m z_i^*L_i$, computed using the dual variable $z^*$ corresponding to
the reference optimal solution.
AdGDC$(a,b)$ sometimes takes stepsizes exceeding these values.

\begin{figure}[H]
 \centering
 \includegraphics[width=.49\linewidth]
  {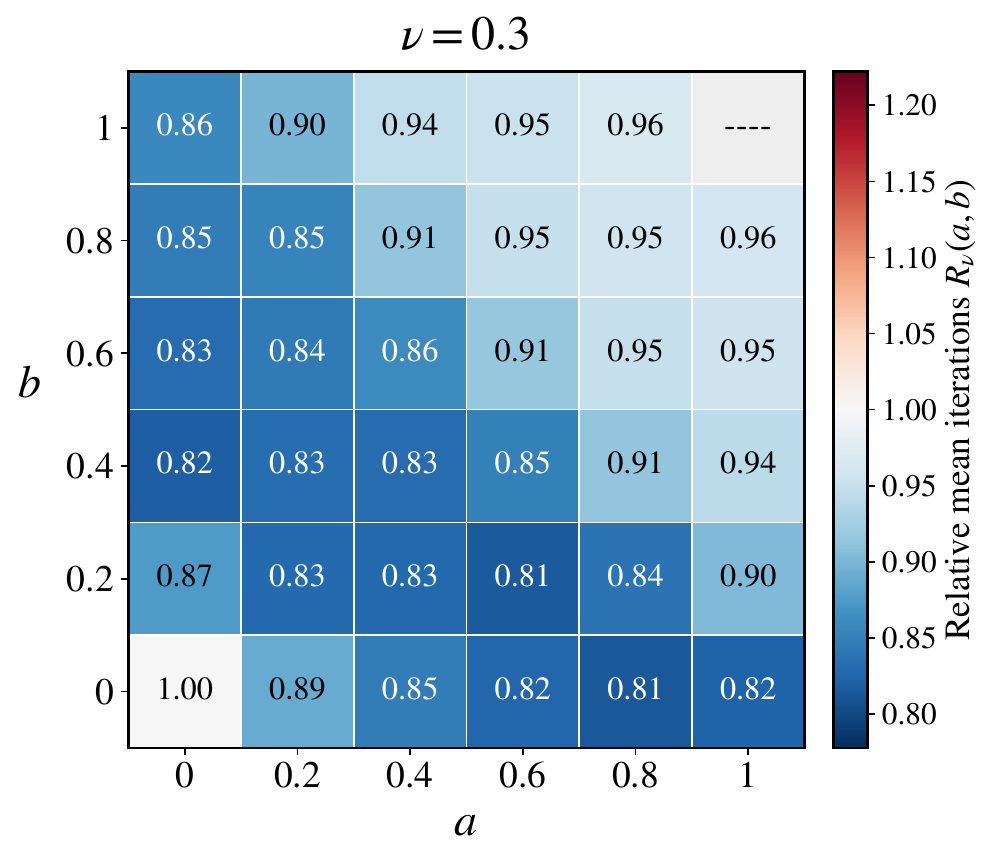}
 \hfill
 \includegraphics[width=.49\linewidth]
  {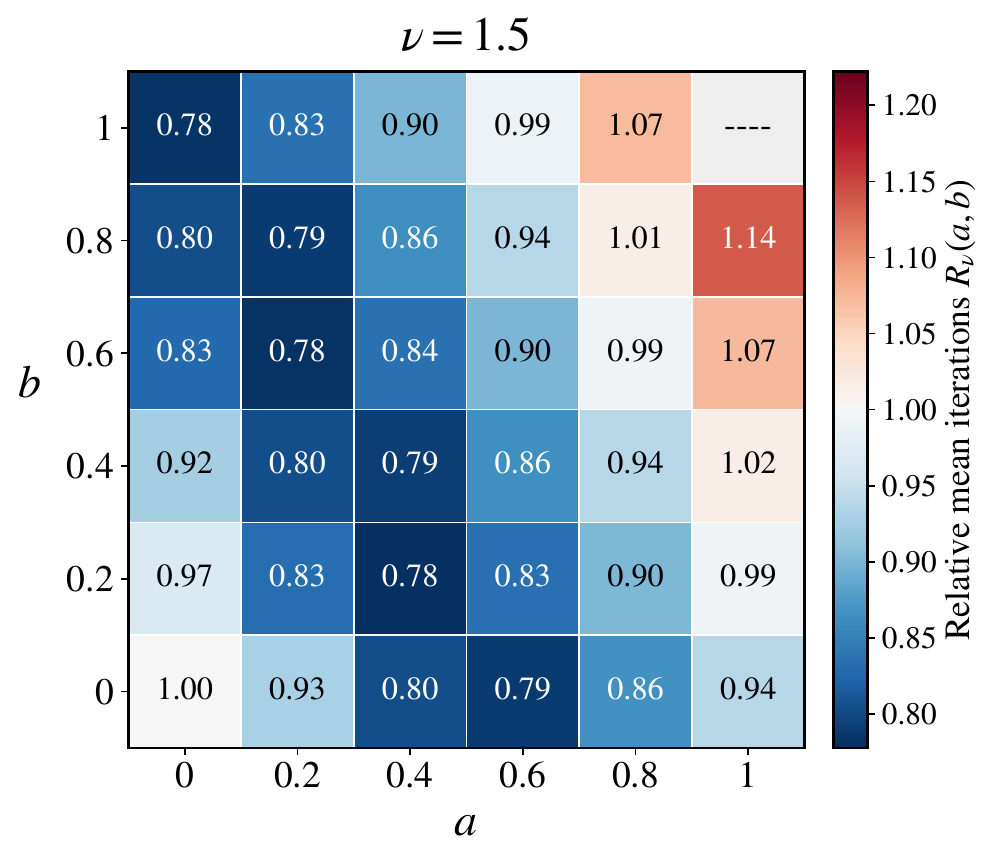}
\caption{Grid search over the correction coefficients $(a,b)$.
Each cell gives the relative mean iteration count $R_\nu(a,b)$ with respect to AdGDC$(0,0)$.
Left: $\nu=0.3$.
Right: $\nu=1.5$.}
 \label{fig:exp_grid}
\end{figure}

Next, to examine sensitivity to the correction coefficients $(a,b)$, we performed a grid search
over the 35 pairs
\[
 (a,b)\in
 \{0,0.2,0.4,0.6,0.8,1\}^2\cap\Omega,
\]
where $\Omega$ is the parameter set defined in~\eqref{eq:def_Omega} and $(1,1)$ is excluded.
For each pair we used the same problem instances and starting points as above and computed the
mean iteration count over 100 runs.
All runs in the grid search satisfied the stopping criterion.
Using the mean iteration count of AdGDC$(0,0)$ as the baseline for each $\nu$, we computed
\[
 R_\nu(a,b):=
 \dfrac{\overline \#_\nu(a,b)}
 {\overline \#_\nu(0,0)},
\]
where $\overline \#_\nu(a,b)$ is the mean iteration count of AdGDC$(a,b)$ over 100 runs.
Thus $R_\nu(a,b)<1$ indicates fewer iterations on average than AdGDC$(0,0)$.
Figure~\ref{fig:exp_grid} shows the results.

At $\nu=0.3$, the smallest mean iteration count was 11.34, achieved by $(a,b)=(0.8,0)$ with
relative value $0.811$. 
At $\nu=1.5$, $(a,b)=(0.4,0.2)$ gave the smallest mean count, 10.47, with relative value $0.778$. 
In both settings several parameter choices reduced the mean count relative to AdGDC$(0,0)$. 
However, large values of both correction coefficients did not consistently yield good performance.

\subsection{Discussion}
\label{subsec:experiments_discussion}

In these experiments the correction terms in AdGDC$(a,b)$ sometimes reduced the mean iteration
count relative to the baseline AdGDC$(0,0)$. 
Which of the two terms was more effective depended on the problem setting. 
In particular, among the four representative choices compared in
Table~\ref{tab:exp_method_summary}, AdGDC$(1,0)$ had the smallest mean count at $\nu=0.3$, whereas
AdGDC$(0,1)$ had the smallest at $\nu=1.5$.
As shown in Section~\ref{sec:parameter_variants},
$\ip{y^k}{s^k}$ and $\sigma_k$ are both nonnegative correction terms, but neither dominates the
other in general.
The relative performance of the methods using the respective terms also varied with the problem
setting in our experiments, suggesting that the terms can play different roles.

In Figure~\ref{fig:exp_stepsizes}, AdGDC$(a,b)$ sometimes takes stepsizes exceeding $1/L$. 
The value $1/L$ is determined from upper bounds on the Lipschitz constants of the gradients of the
functions $F_i$. 
It serves as a common empirical benchmark, but does not impose an upper bound on the stepsizes of AdGDC$(a,b)$. 
Using information obtained during the iterations allows our method to take steps larger than a
fixed stepsize based on a global Lipschitz constant. 
The smaller iteration counts of AdGDC$(a,b)$ relative to the fixed-step prox-linear and
accelerated prox-linear methods in these experiments are consistent with the effectiveness of
adaptive stepsize selection. 

In the grid search in Figure~\ref{fig:exp_grid}, the iteration count did not decrease
monotonically as the correction coefficients $(a,b)$ increased. 
For example, at $\nu=1.5$, AdGDC$(0,0)$ required 13.46 iterations on average, while AdGDC$(1,0.8)$
required 15.29, more than the method without correction. 
AdGDC$(0.4,0.2)$, on the other hand, required 10.47. 
Thus assigning large weights to both corrections need not accelerate convergence. 
This observation does not contradict Proposition~\ref{prop:parameter_monotonicity}, which states
that, for fixed iteration data, increasing $a,b$ cannot decrease the stepsize upper bound
$\widehat t_k$. 
The proposition does not imply any monotonicity in the iterates generated with different parameter
choices or in the number of iterations required for termination.
Even if a larger correction permits a larger stepsize, the resulting change in the iterates may
slow convergence.
Enlarging the admissible stepsize range and improving computational efficiency are therefore
distinct considerations.

In these experiments, using one or both correction terms in AdGDC$(a,b)$ could reduce the
iteration count to convergence, and both the choice of correction and its weight affected the
convergence speed. 
These observations concern the present experimental setting.
The relative performance may differ for other problem classes. 
Likewise, the best parameter choices identified by the grid search in Figure~\ref{fig:exp_grid}
pertain to the instances tested. 
The relation between problem structure and suitable correction coefficients warrants further
investigation.

\section{Conclusions}
\label{sec:conclusion}

We proposed AdGDC$(a,b)$, a prox-linear method with a linesearch-free adaptive stepsize for convex
composite problems involving a nonlinear inner mapping and a possibly nonsmooth convex outer
function. 
Building on the adaptive gradient method of Malitsky and Mishchenko, the method combines two
nonnegative quantities arising from the composite structure through the correction coefficients
$(a,b)$. 
This provides a common framework for stepsize rules based on different correction terms. 
The method determines the stepsize from iteration data without a prior Lipschitz estimate or a
line search.

We derived a fundamental inequality common to the AdGDC$(a,b)$ family and proved boundedness of
the iterates under basic convexity assumptions. 
For AdGDC$(a,b)$, which uses two nonnegative quantities to relax the stepsize condition, we
established, throughout $(a,b)\in\Omega$, a uniform positive lower bound on the stepsizes, an
$O(1/K)$ best-iterate convergence rate, convergence of the objective values, and convergence of
the entire sequence.

In our numerical experiments we compared AdGDC$(a,b)$ with the prox-linear and accelerated
prox-linear methods~\cite{DrusvyatskiyPaquette2019} on a nonsmooth convex composite problem
defined by the maximum of smooth convex functions based on the log-cosh loss. 
The proposed method performed favorably against the existing methods.
For some choices of $(a,b)\ne(0,0)$, AdGDC$(a,b)$ required fewer iterations
than the baseline AdGDC$(0,0)$.
The relative effectiveness of the two corrections depended on the problem setting.
Increasing the coefficients $(a,b)$ did not consistently accelerate convergence, suggesting that
an appropriate correction strength also depends on the problem characteristics.

Future work includes evaluating performance across different problem structures and scales,
identifying the problem characteristics that determine the effectiveness of the two corrections,
and developing adaptive choices of $(a,b)$. 
Extending the method to more general composite optimization problems is another important
direction.

\appendix

\section{Proofs of convergence of the objective values and the entire sequence}
\label{app:convergence_only}

We prove Lemma~\ref{lem:objective_convergence} and Theorem~\ref{thm:global}.
Throughout, fix $(a,b)\in\Omega$, suppose that AdGDC$(a,b)$ does not terminate after finitely many
iterations, and assume that (L) holds.

By Theorem~\ref{thm:bounded}, the iterates $\{x^k\}$ and $\{s^k\}$ are bounded.
By Theorem~\ref{thm:positive_step}, there exists $t_{\inf}>0$ such that $t_k\ge t_{\inf}$ for
every $k\ge0$.
Moreover, the estimate in Corollary~\ref{cor:summability} for the residual $e^k$ defined
in~\eqref{eq:ek} and the estimate in Lemma~\ref{lem:y_weighted_summability} for $y^k$ give
\begin{align}
 \sum_{k=1}^{\infty}\norm{e^k}^2&<\infty,
 \label{eq:joint_core_sums}\\
 \sum_{k=1}^{\infty}t_k\theta_k\norm{y^k}^2&<\infty .
 \label{eq:joint_weighted_y}
\end{align}
We use these results below.

\subsection{Convergence of the objective values}

\begin{proof}[Proof of Lemma~\ref{lem:objective_convergence}]
Fix $x^*\in X^* $.
For an arbitrary $\varepsilon>0$, suppose that the set
\[
 \mathcal K_\varepsilon
 :=\{k\ge1:\delta_k\ge\varepsilon\}
\]
is infinite.

In Theorem~\ref{thm:bounded},
take $K=k$ in~\eqref{eq:telescoped_rate}
and discard nonnegative terms to obtain
\[
 2(\omega_k+t_{k+1}\theta_{k+1})\delta_k
 \le D_{x^*}.
\]
Here $D_{x^*}$ is the constant defined in~\eqref{eq:D}.
The definition of $\omega_k$ in \eqref{eq:omega} gives
\[
 2t_k(1+\theta_k)\delta_k\le D_{x^*}.
\]
Thus, on $k\in\mathcal K_\varepsilon$, it follows from $\theta_k>0$ that
\begin{equation}\label{eq:joint_t_eps}
 t_k\le\frac{D_{x^*}}{2\varepsilon}
 =:T_\varepsilon.
\end{equation}

Next consider $v^k=s^{k+1}/t_k$ defined in~\eqref{eq:vk}.
For the model $m_x$ defined in~\eqref{eq:model},
the optimality condition~\eqref{eq:optim_condition} for subproblem~\eqref{eq:subproblem} gives
$-v^k\in\partial m_{x^k}(s^{k+1})$.
Thus convexity of $m_{x^k}$ and~
\eqref{eq:model_minorant} imply
\[
  h^* 
 \ge m_{x^k}(x^*-x^k)
 \ge g(\widetilde F^k)
     +\ip{v^k}{x^{k+1}-x^*}.
\]
Hence
\begin{align*}
 \delta_k
 &\le h(x^k)-g(\widetilde F^k)
       -\ip{v^k}{x^{k+1}-x^*}\\
 &=g(F(x^k))
   -g(F(x^k)+\J(x^k)^\top s^{k+1})
   -\ip{v^k}{x^{k+1}-x^*}.
\end{align*}

Local Lipschitz continuity of $g$, boundedness of the iterates, and boundedness of $\J(x^k)$ imply
that there exist constants $A,B>0$ such that
\[
 \delta_k
 \le A\norm{s^{k+1}}+B\norm{v^k}
 =(At_k+B)\norm{v^k}.
\]
Thus, on $k\in\mathcal K_\varepsilon$,
it follows from~\eqref{eq:joint_t_eps} that
\begin{equation}\label{eq:joint_velocity_gap}
 \norm{v^k}
 \ge\frac{\varepsilon}{AT_\varepsilon+B}
 =:\xi_\varepsilon>0.
\end{equation}

Fix an arbitrary positive integer $N$.
By Lemma~\ref{lem:theta_bound}, we have $\theta_k\le\varphi:=(1+\sqrt5)/2$, so for
$k\in\mathcal K_\varepsilon$,
\[
 \begin{aligned}
 t_{k+j}
 &\le\varphi^jT_\varepsilon
 &&(0\le j\le N),\\
 t_{k+j}\theta_{k+j}
 &=\frac{t_{k+j}^2}{t_{k+j-1}}
 \ge\frac{t_{\inf}^2}
          {\varphi^{j-1}T_\varepsilon}
 &&(1\le j\le N).
 \end{aligned}
\]

Hence, by~\eqref{eq:joint_weighted_y}, for each fixed $j=1,\ldots,N$ we have
\[
 y^{k+j}\longrightarrow0
 \qquad
 (k\to\infty,\ k\in\mathcal K_\varepsilon).
\]
Moreover, \eqref{eq:joint_core_sums} gives $e^k\to0$, and $t_k\ge t_{\inf}>0$ implies
$e^k/t_k\to0$.

By Lemma~\ref{lem:velocity_difference},
a telescoping sum of~\eqref{eq:velocity_difference} gives
\[
 v^{k+j}-v^k
 =-\sum_{\ell=1}^j y^{k+\ell}
 +\sum_{\ell=1}^j\frac{e^{k+\ell}}{t_{k+\ell}}.
\]
Since $N$ is fixed,
\begin{equation}\label{eq:joint_block_velocity}
 \max_{0\le j\le N}\norm{v^{k+j}-v^k}
 \longrightarrow0
 \qquad
 (k\to\infty,\ k\in\mathcal K_\varepsilon).
\end{equation}

By boundedness of $\{x^k\}$, choose $M>0$ such that $\norm{x^k}\le M$ for all $k\ge0$.
Fix a positive integer $N$ satisfying
\[
 N>\frac{4M}{t_{\inf}\xi_\varepsilon}.
\]
For $k\in\mathcal K_\varepsilon$,
since~\eqref{eq:joint_velocity_gap} gives $v^k\ne0$, we may define
\[
 w^k:=\frac{v^k}{\norm{v^k}}.
\]
By~\eqref{eq:joint_block_velocity},
for sufficiently large $k\in\mathcal K_\varepsilon$,
\[
 \begin{aligned}
 \ip{v^{k+j}}{w^k}
 &=\norm{v^k}
   +\ip{v^{k+j}-v^k}{w^k}\\
 &\ge\xi_\varepsilon
       -\norm{v^{k+j}-v^k}\\
 &\ge\frac{\xi_\varepsilon}{2}
 \qquad(0\le j\le N-1).
 \end{aligned}
\]

Combining this with $t_k v^k = s^{k+1}=x^{k+1}-x^k$ yields
\[
 \begin{aligned}
 \ip{x^{k+N}-x^k}{w^k}
 &=\sum_{j=0}^{N-1}
   t_{k+j}\ip{v^{k+j}}{w^k}\\
 &\ge\frac{Nt_{\inf}\xi_\varepsilon}{2}\\
 &>2M.
 \end{aligned}
\]
However, $\norm{w^k}=1$ and $\norm{x^k}\le M$ give
\[
 \ip{x^{k+N}-x^k}{w^k}
 \le\norm{x^{k+N}-x^k}
 \le2M,
\]
a contradiction.

Thus $\mathcal K_\varepsilon$ is finite.
Since $\varepsilon>0$ was arbitrary, $\delta_k\to0$ holds.
\end{proof}

\subsection{Convergence of the entire sequence}

\begin{proof}[Proof of Theorem~\ref{thm:global}]
By Lemma~\ref{lem:objective_convergence}, we have $\delta_k\to0$.
Since the iterates $\{x^k\}$ are bounded and $h$ is continuous, $\{x^k\}$ has at least one cluster
point, and all its cluster points belong to $ X^* $.

For any $x^*\in X^* $, define
\[
 \mathcal V_k(x^*)
 :=\norm{x^k-x^*}^2+\norm{s^k}^2
   +2t_k\theta_k\delta_{k-1}
 \qquad(k\ge1).
\]

By~\eqref{eq:strong_onestep},
\[
 \mathcal V_{k+1}(x^*)
 +2\omega_k\delta_k
 +2\norm{e^k}^2
 +2t_k\theta_k
 \left(
 (1-a)\ip{y^k}{s^k}
 +(1-b)\sigma_k
 \right)
 \le\mathcal V_k(x^*).
\]

Here~\eqref{eq:omega} gives $\omega_k\ge0$, \eqref{eq:delta_def} gives $\delta_k\ge0$,
\eqref{eq:ymonotone} gives $\ip{y^k}{s^k}\ge0$, \eqref{eq:prev_gap} gives $\sigma_k\ge0$,
and~\eqref{eq:def_Omega} gives
$1-a\ge0$ and $1-b\ge0$.
Hence
\[
 0\le
 \mathcal V_{k+1}(x^*)
 \le\mathcal V_k(x^*).
\]
Thus $\{\mathcal V_k(x^*)\}$ is nonnegative and nonincreasing and therefore has a limit for every
$x^*\in X^* $.

Consider two cluster points $x_{\mathrm{A}}^*,x_{\mathrm{B}}^*\in X^* $ of $\{x^k\}$.
We have
\[
 \begin{aligned}
 \mathcal V_k(x_{\mathrm{A}}^*)-\mathcal V_k(x_{\mathrm{B}}^*)
 &=\norm{x^k-x_{\mathrm{A}}^*}^2-\norm{x^k-x_{\mathrm{B}}^*}^2
 \end{aligned},
\]
and the left-hand side has a limit as $k\to\infty$. 
Along a subsequence converging to $x_{\mathrm{A}}^*$, the right-hand side converges to
$-\norm{x_{\mathrm{A}}^*-x_{\mathrm{B}}^*}^2$, whereas along a subsequence converging to
$x_{\mathrm{B}}^*$ it converges to $\norm{x_{\mathrm{A}}^*-x_{\mathrm{B}}^*}^2$. 
By uniqueness of the limit on the left, we have
\[
 -\norm{x_{\mathrm{A}}^*-x_{\mathrm{B}}^*}^2=\norm{x_{\mathrm{A}}^*-x_{\mathrm{B}}^*}^2,
\]
which gives $x_{\mathrm{A}}^*=x_{\mathrm{B}}^*$.

Thus the cluster point $x^*\in X^* $ is unique, and
\[
 x^k\longrightarrow x^*
 \qquad(k\to\infty).
\]
In particular, $s^k=x^k-x^{k-1}\to0$.
\end{proof}

\section{Existing methods for composite problems}
\label{app:prox_linear}

This appendix describes the two comparison methods used in the numerical experiments in
Section~\ref{sec:experiments}: the prox-linear method~\cite[Algorithm~1]{DrusvyatskiyPaquette2019}
and the accelerated prox-linear method~\cite[Algorithm~7]{DrusvyatskiyPaquette2019} of
Drusvyatskiy and Paquette. 
We describe the algorithms for problem~\eqref{eq:exp_penalty_objective} and then detail their
implementation in our experiments. 
The descriptions below follow those in~\cite{DrusvyatskiyPaquette2019}, with the notation adapted
as necessary for the present paper.

For a point $x\in\R^n$ and stepsize $t>0$, define
\begin{equation}\label{eq:baseline_direction}
 d_t(x):=
 \argmin_{d\in\R^n}
 \left\{
 m_x(d)+\frac{1}{2t}\norm d^2
 \right\},
\end{equation}
where $m_x$ is the linearized model defined in~\eqref{eq:model}.
Subproblem~\eqref{eq:baseline_direction} is strongly convex and has a unique solution $d_t(x)$.

The prox-linear method of Drusvyatskiy and Paquette~\cite[Algorithm~1]{DrusvyatskiyPaquette2019}
uses a fixed stepsize $t>0$ and updates the iterate by
\[
 x^{k+1}=x^k+d_t(x^k).
\]
In our experiments, we use $L=\max_iL_i$ defined in~\eqref{eq:exp_L} and set $t=1/L$.
Since $L_i$ bounds the Lipschitz constant of the gradient of each $F_i$, for any $x,d\in\R^n$ we
have
\[
 F_i(x+d)\le
 F_i(x)+\ip{\nabla F_i(x)}{d}
 +\frac{L}{2}\norm d^2.
\]
Thus, for the present problem with $g(u)=\max_i u_i$,
\[
 h(x+d)\le
 m_x(d)+\frac{L}{2}\norm d^2.
\]
In other words, with $t=1/L$, the objective in subproblem~\eqref{eq:baseline_direction} majorizes
the original objective.
We use the subproblem solution directly as the update direction without a line search.

\renewcommand{\thealgorithm}{B.\arabic{algorithm}}
\setcounter{algorithm}{0}

\begin{algorithm}[H]
\caption{Prox-linear method~\cite[Algorithm~1]{DrusvyatskiyPaquette2019}}
\label{alg:baseline_prox_linear}
\begin{algorithmic}[1]
\Require Initial point $x^0\in\R^n$, fixed stepsize $t=1/L$
 \For{$k=0,1,2,\ldots$}
  \State $\displaystyle
   d^k\gets\argmin_{d\in\R^n}
   \left\{
   m_{x^k}(d)+\frac{1}{2t}\norm d^2
   \right\}$
  \State $x^{k+1}\gets x^k+d^k$
 \EndFor
\end{algorithmic}
\end{algorithm}

We next describe the accelerated prox-linear method of Drusvyatskiy and
Paquette~\cite[Algorithm~7]{DrusvyatskiyPaquette2019}.
It maintains two sequences, $\{x^k\}$ and $\{q^k\}$, and solves a prox-linear subproblem based at
an interpolated point.
Fix the initial point $x^0,q^0\in\R^n$ and the parameter $\widetilde\mu>L$.
For $k\ge1$, set
\begin{align*}
 \alpha_k&:=\frac{2}{k+1},\\
 \widehat x^k&:=\alpha_kq^{k-1}+
 (1-\alpha_k)x^{k-1}.
\end{align*}
Update the sequence $\{x^k\}$ by
\begin{equation}\label{eq:baseline_acc_x}
 x^k:=
 \argmin_{x\in\R^n}
 \left\{
 m_{\widehat x^k}(x- \widehat x^k)
 +\frac{\widetilde\mu}{2}\norm{x- \widehat x^k}^2
 \right\},
\end{equation}
and update $\{q^k\}$ by solving the following subproblem, in which $\alpha_k$ modifies the
linearized model:
\begin{equation}\label{eq:baseline_acc_v}
 q^k:=\argmin_{q\in\R^n}
 \Biggl\{
 \frac{1}{\alpha_k}
 g\Bigl(
 F(\widehat x^k)+\alpha_k
 \J(\widehat x^k)^\top(q-q^{k-1})
 \Bigr)+\frac{\widetilde\mu\alpha_k}{2}
 \norm{q-q^{k-1}}^2
 \Biggr\}.
\end{equation}
These are the updates of Drusvyatskiy and Paquette~\cite[Algorithm~7]{DrusvyatskiyPaquette2019}
specialized to problem~\eqref{eq:problem_intro}.

In general this algorithm requires solving both subproblems~\eqref{eq:baseline_acc_x}
and~\eqref{eq:baseline_acc_v}.
For our problem they reduce to the same search-direction subproblem.
Indeed, in~\eqref{eq:baseline_acc_v}, set
\[
 d=\alpha_k(q-q^{k-1}).
\]
The resulting objective is
\[
 \frac{1}{\alpha_k}
 \left\{
 g\bigl(F(\widehat x^k)+\J(\widehat x^k)^\top d\bigr)
 +\frac{\widetilde\mu}{2}\norm d^2
 \right\}.
\]
Since $\alpha_k>0$, minimizing this objective is equivalent to solving
\begin{equation}\label{eq:baseline_acc_common}
 d^k:=
 \argmin_{d\in\R^n}
 \left\{
 m_{\widehat x^k}(d)
 +\frac{\widetilde\mu}{2}\norm d^2
 \right\}.
\end{equation}
The change of variables $d=x- \widehat x^k$ in~\eqref{eq:baseline_acc_x} gives the same
subproblem~\eqref{eq:baseline_acc_common}.
We may therefore update the iterates by
\[
 x^k=\widehat x^k+d^k,\qquad
 q^k=q^{k-1}+\frac{d^k}{\alpha_k}.
\]
Thus, in our setting, the accelerated method requires only one subproblem solve per iteration.

In our experiments we set $q^0=x^0$ and use $\widetilde\mu=1.01L$.
The accelerated prox-linear method used in the experiments is summarized below.

\begin{algorithm}[H]
\caption{Accelerated prox-linear method~\cite[Algorithm~7]{DrusvyatskiyPaquette2019}}
\label{alg:baseline_acc_prox_linear}
\begin{algorithmic}[1]
\Require Initial point $x^0\in\R^n$,
$q^0\gets x^0$, parameter $\widetilde\mu=1.01L$
 \For{$k=1,2,\ldots$}
  \State $\alpha_k\gets2/(k+1)$
  \State $\widehat x^k\gets
   \alpha_kq^{k-1}+(1-\alpha_k)x^{k-1}$
  \State $\displaystyle
   d^k\gets\argmin_{d\in\R^n}
   \left\{
   m_{\widehat x^k}(d)+\frac{\widetilde\mu}{2}\norm d^2
   \right\}$
  \State $x^k\gets \widehat x^k+d^k$
  \State $q^k\gets q^{k-1}+d^k/\alpha_k$
 \EndFor
\end{algorithmic}
\end{algorithm}

Finally, we describe how the subproblems are solved in our experiments.
Since $g(u)=\max_{1\le i\le m}u_i$, subproblem~\eqref{eq:baseline_direction} can be written as
\[
 \min_{d\in\R^n}
 \left\{
 \max_{1\le i\le m}
 \bigl(F_i(x)+\ip{\nabla F_i(x)}{d}\bigr)
 +\frac{1}{2t}\norm d^2
 \right\}.
\]
Using
$\Delta_m:=\{z\in\R_+^m:\sum_{i=1}^m z_i=1\}$
gives the dual problem
\[
 \min_{z\in\Delta_m}
 \left\{
 \frac{t}{2}\norm{\J(x)z}^2
 -\ip{F(x)}{z}
 \right\}.
\]
If $z_t(x)$ is a dual optimal solution, the primal optimal solution is given by
\[
 d_t(x)=-t\J(x)z_t(x).
\]
In the experiments we solve this convex quadratic program with OSQP and update the iterate using
the resulting direction.
For the accelerated prox-linear method, \eqref{eq:baseline_acc_common} is an instance of
subproblem~\eqref{eq:baseline_direction} with $t=1/\widetilde\mu$, so the same dual formulation
and QP solver apply.
Both methods use the stopping criterion~\eqref{eq:exp_stop}.

\end{document}